\pdfoutput=1
\documentclass[11pt]{article}
\usepackage[letterpaper,margin=1in]{geometry}
\usepackage{graphicx}
\usepackage[utf8]{inputenc}
\usepackage[T1]{fontenc}
\usepackage{mathptmx}
\usepackage{amsmath}
\usepackage{amssymb}
\usepackage{float}
\usepackage[tableposition=top]{caption}
\usepackage[numbers,square,comma,sort&compress]{natbib}
\usepackage[colorlinks=true,linkcolor=blue,citecolor=blue,urlcolor=blue]{hyperref}

\makeatletter
\newenvironment{ruledtabular}{%
  \def\tabular##1{\csname tabular*\endcsname{\linewidth}{@{\extracolsep{\fill}}##1}\hline\hline}%
  \def\endtabular{\hline\hline\csname endtabular*\endcsname}%
}{}
\newenvironment{acknowledgments}{\section*{Acknowledgments}}{}
\g@addto@macro\appendix{%
  \@addtoreset{equation}{section}%
}
\AtBeginDocument{\immediate\write\@auxout{\string\citation{REVTEX41Control}}\immediate\write\@auxout{\string\citation{aip41Control}}}
\makeatother

\begin{document}

\title{Causal Fingerprints Reveal Hidden Mechanisms in Partially-Observed Dynamical Systems}

\author{Zhang Jiang, Marios Andreou, and Nan Chen\thanks{Corresponding author: \href{mailto:chennan@math.wisc.edu}{chennan@math.wisc.edu}}\\[4pt]
\normalsize Department of Mathematics, University of Wisconsin--Madison, Madison, WI 53706, USA}

\date{}

\maketitle

\begin{abstract}
In partially-observed dynamical systems, distinct hidden mechanisms can produce similar observable summary statistics while explaining individual events differently. Matching statistics therefore does not guarantee matching mechanisms. We introduce causal fingerprints: realization-dependent patterns that track how information about a hidden state is revealed through subsequent evolution. Within a prescribed dynamical model, we compare the filtering posterior, based on observations up to a given time, with the smoothing posterior, obtained after adding an expanding window of future observations. We quantify this retrospective information gain using Kullback-Leibler divergence and track its evolution as the future window of observational data expands. The fingerprint reveals both the strength and temporal reach of hidden influence, which can vary across events and their stages. Because the observed trajectory is held fixed while hidden-state inference depends on the candidate model, fingerprints reveal how competing dynamical mechanisms explain the same observations even when selected statistics agree. In linear Gaussian systems, we derive a characteristic fingerprint timescale for the temporal reach of hidden influence. Calibrated surrogates matching the observable mean, variance, and decorrelation time can still produce distinct fingerprints when feedback memory, hidden-noise memory, or hidden forcing timescales are misrepresented. Nonlinear examples reveal mechanism-dependent signatures for extreme bursts, oscillatory variability, and regime switching, including multiplicative instability versus additive forcing, deterministic nonlinear feedback versus state-dependent multiplicative noise, and intrinsic bistability versus hidden double-well forcing. Causal fingerprints therefore provide a time- and event-resolved diagnostic for asking not only whether a model reproduces observed statistics, but how its hidden dynamics explain individual events.
\end{abstract}

\begin{quotation}
Models with distinct mechanisms can reproduce similar statistical patterns yet explain the same observed event differently. Telling these mechanistic explanations apart is particularly difficult when the variables that drive an event are hidden from observation. However, a hidden variable that drives the later evolution of the observed variables leaves information in subsequent observations, so these observations can refine the inference about the earlier state of the hidden variable. We quantify this information gain as causal fingerprints, which reveal how strongly hidden variables influence the observed evolution and how long this influence persists across the stages of individual events. In examples involving extreme bursts, oscillatory variability, and regime switching, models with closely matched statistics produce distinct fingerprints along the same observed trajectories, showing that fingerprints reveal the different mechanisms through which these models explain each event.
\end{quotation}

\section{Introduction}

Complex dynamical systems arise throughout science and engineering, including geophysical dynamics~\cite{carrassi2018data}, turbulence~\cite{majdaharlim2012filtering}, neuroscience~\cite{paninski2010new}, epidemiology~\cite{king2016statistical}, and financial markets~\cite{heston1993closed}. A defining feature of many such systems is that they are {partially observed}: in practice, the observed variables represent only a subset of the full dynamical state, while many relevant variables remain hidden or unresolved due to limited observations, coarse-grained modeling, or intrinsic spatiotemporal multiscale structures~\cite{chenmajda2018conditional}. These hidden variables often represent slow environmental drivers, unresolved sub-grid or small-scale processes~\cite{majdaharlim2012filtering}, or latent regulatory mechanisms~\cite{paninski2010new}.

A central challenge in studying partially-observed systems is that different hidden mechanisms can produce very similar observable statistics, such as probability distributions, autocorrelations, or energy spectra~\cite{majdachen2018modelerror, majda2018strategies}. As a result, model identification based solely on observable statistics can admit different explanations for the underlying latent dynamics~\cite{villaverde2016structural}. Understanding the causal influence that hidden variables have on the observables' evolution is therefore essential, providing insight into the mechanisms that drive observable variability and shape phenomena such as extreme events~\cite{andreou2026assimilative}, oscillatory variability~\cite{vialard2025recharge}, and regime transitions~\cite{hanggi1990reaction}. In dynamical systems, such influence can vary substantially along a trajectory and across different stages of an event, making time-local and event-resolved diagnosis particularly important.

Existing causal inference frameworks have been successfully applied to a wide range of dynamical systems. Methods such as Granger causality~\cite{granger1969investigating}, transfer entropy~\cite{schreiber2000measuring}, and graphical-model-based approaches~\cite{runge2019inferring} are particularly useful for identifying statistical or predictive relationships among observed, or reliably reconstructed, variables. However, their considered setting poses two significant difficulties. First, some of the variables responsible for the observed evolution are themselves hidden. Second, their influence may change substantially along a trajectory, so relationships inferred from long records can miss mechanisms that become important only during a particular event or stage of an event. Our interest is therefore not only in whether a causal relationship exists on average, but in how hidden influence unfolds during an individual realization.

A complementary line of work uses dynamical structures more explicitly. Information-flow formulations \cite{liang2013liang}, perturbation experiments \cite{meinshausen2016methods}, and response-theory approaches \cite{falasca2025probing}, including fluctuation-dissipation methods, use the governing dynamics to characterize directional entropy transfer, sensitivity, or response to perturbations. The question considered here is different: an observed trajectory is given while important causal variables remain hidden, with candidate models providing competing explanations for the latent dynamics behind that trajectory.

This motivates the concept of {causal fingerprints} for partially-observed dynamical systems. Distinct latent mechanisms can produce similar observable statistics, making them difficult to distinguish using probability distributions, long-time correlations, or related diagnostics~\cite{schreiber2000surrogate,majdachen2018modelerror}. Yet they need not explain an individual trajectory in the same way. Hidden influence may become important at different stages of an event or remain relevant over different future horizons. We refer to these realization-dependent temporal patterns of hidden-to-observed influence as causal fingerprints.

Building on assimilative causal inference and related work on causal influence ranges~\cite{andreou2026assimilative,andreou2026bridging}, we use these realization-dependent fingerprints to diagnose competing dynamical mechanisms. Fingerprints computed under different candidate models can be compared along the same observed trajectory, revealing differences in the hidden dynamical explanations assigned to the observations.

To quantify these patterns, we adopt a {Bayesian inverse-problem perspective of causality}~\cite{andreou2026assimilative}. When a hidden state affects the subsequent observed evolution, later observations can sharpen the inference of that earlier state. Using {data assimilation} \cite{law2015data, sanz2023inverse}, we compare the filtering posterior based on observations up to the current time with the smoothing posterior obtained after incorporating an expanding window of future observations~\cite{sarkka2013bayesian}. Extending this future window tracks how far information about the hidden state persists in the subsequent observations. Repeating the comparison along the trajectory shows how this pattern changes across events and during their various stages of onset, development, peak, and decay. The diagnostic is therefore observation-conditioned but model-informed: the observed trajectory is fixed, while hidden-state inference depends on the candidate dynamical model. Fingerprints computed from the same observations can consequently reveal how different hidden mechanisms explain the same observed evolution.

The remainder of the paper is organized as follows. Section~\ref{sec:theory} presents the mathematical formulation, definition, and interpretation of causal fingerprints. Section~\ref{sec:linear} develops the analyses and mechanism tests for linear Gaussian systems. Section~\ref{sec:nonlinear} examines nonlinear dynamical systems involving extreme events, oscillatory variability, and regime switching. Section~\ref{sec:discussion} summarizes the numerical findings and discusses future directions. Additional derivations and exposition can be found in the Appendix.

\section{Causal Fingerprints in Partially-Observed Systems}
\label{sec:theory}

This section develops the causal-fingerprint framework for partially-observed dynamical systems. The central idea is to compare what can be inferred about a hidden state from observations available up to a given time with what can be inferred about the same state after subsequent observations are added. If that hidden state affects the later observed evolution, those later observations should carry information about it.

\subsection{Partially-observed dynamical systems}

We consider dynamical systems in which only a subset of the state variables is observed, while the remaining variables are hidden. Such systems naturally arise when observations are limited or when unresolved processes influence the observed dynamics, as in numerical weather forecasting~\cite{carrassi2018data}, feedback control~\cite{astrommurray2008feedback}, and turbulent flow modeling~\cite{majdaharlim2012filtering}.

Let $\mathbf{X}_t:=\mathbf{X}(t) \in \mathbb{R}^{d_X}$ denote the observed variables and $\mathbf{Y}_t:=\mathbf{Y}(t) \in \mathbb{R}^{d_Y}$ denote the hidden variables at state time \(t\geq 0\). Without a time subscript, \(\mathbf{X}\) and \(\mathbf{Y}\) denote the corresponding processes or trajectories. We consider the stochastic dynamical system
\begin{equation}
\begin{aligned}
d\mathbf{X}_t &= \mathbf{f}(t,\mathbf{X}_t,\mathbf{Y}_t)\,dt + \boldsymbol{\sigma}_X(t,\mathbf{X}_t,\mathbf{Y}_t)\,d\mathbf{W}_X,\\
d\mathbf{Y}_t &= \mathbf{g}(t,\mathbf{X}_t,\mathbf{Y}_t)\,dt + \boldsymbol{\sigma}_Y(t,\mathbf{X}_t,\mathbf{Y}_t)\,d\mathbf{W}_Y .
\end{aligned}
\label{eq:system}
\end{equation}
where $\mathbf{f}$ and $\mathbf{g}$ are generally nonlinear drift functions, $\boldsymbol{\sigma}_X$ and $\boldsymbol{\sigma}_Y$ are possibly state-dependent diffusion coefficients, and $\mathbf{W}_X$ and $\mathbf{W}_Y$ are independent Wiener processes. Throughout the paper, all SDEs are interpreted in the It\^o sense~\cite{sarkka2019applied}.

\subsection{Definition of the causal fingerprint}

Consider a prescribed dynamical model of the form \eqref{eq:system} (hereafter called the inference model), which we use to infer the hidden state from an observed trajectory. At state time $t$, the filtering distribution
$p(\mathbf{Y}_t\mid \mathbf{X}_{[0,t]})$
describes what can be inferred about $\mathbf{Y}_t$ from the observations available up to that time, \mbox{$\mathbf{X}_{[0,t]}=\{\mathbf{X}_s:0\leq s\leq t\}$}. After observations over a future window $(t,t+\tau]$ are added, the corresponding smoothing distribution becomes
$p(\mathbf{Y}_t\mid \mathbf{X}_{[0,t+\tau]})$.
The two conditional distributions concern the same hidden state. Their difference comes only from the information contained in the subsequent observed evolution.

This filter-smoother comparison gives the inverse-problem interpretation of hidden causal influence used here. If $\mathbf{Y}_t$ affects the subsequent evolution of $\mathbf{X}$ under the prescribed dynamics, then the consequences of that influence can leave information in the later observations $\mathbf{X}_{[0,t+\tau]}$. Those observations can therefore revise the inference of the earlier hidden state. The causal interpretation comes from this temporal relation and the dynamical model; the divergence introduced below quantifies the resulting information gain.

For a realized observed trajectory $\mathbf{X}$, we define the causal fingerprint at state time $t$ and future-window length $\tau$ as the smoother-to-filter Kullback-Leibler (KL) divergence \cite{joyce2025kullback, van2014renyi}:
\begin{equation}
\mbox{Causal Fingerprint:}\qquad
\mathcal{I}(t,\tau)
:=
D_{\mathrm{KL}}\left(
p(\mathbf{Y}_t \mid \mathbf{X}_{[0,t+\tau]})
\Vert
p(\mathbf{Y}_t \mid \mathbf{X}_{[0,t]})
\right).
\label{eq:info_gain_natural}
\end{equation}
For probability densities \(p\) and \(q\),
\[
D_{\mathrm{KL}}(p\Vert q)
:=
\int p(\mathbf{y})
\log\frac{p(\mathbf{y})}{q(\mathbf{y})}
\,\mathrm{d}\mathbf{y}.
\]
At $\tau=0$, the filtering and smoothing distributions coincide and $\mathcal I(t,0)=0$.

As the future window length $\tau$ expands, the fingerprint records how long subsequent observations continue to refine the inference of $\mathbf{Y}_t$. It therefore probes the temporal reach of the model-implied influence of the hidden state on the observed evolution~\cite{andreou2026bridging}. The dependence on $t$ makes this diagnostic trajectory-dependent: the instantaneous strength and temporal reach of hidden influence can change from one event to another and across different stages of the same event.

\subsection{Interpretation of the causal fingerprint}

A useful information-theoretic interpretation follows by averaging the fingerprint over observation paths. Under the same inference model, we have
\begin{align}
\mathbb E_{\mathbf{X}_{[0,t+\tau]}}
\bigl[\mathcal I(t,\tau)\bigr]
&=
I\left(
\mathbf{Y}_t;
\mathbf{X}_{(t,t+\tau]}
\mid
\mathbf{X}_{[0,t]}
\right),
\label{eq:expected_fingerprint_cmi}
\end{align}
where the right-hand side is the conditional mutual information~\cite{cover2006elements}. Thus, the trajectory-dependent fingerprint is a pathwise counterpart of a standard conditional information measure. The averaged fingerprint is nondecreasing as the future window expands, whereas $\mathcal I(t,\tau)$ along an individual realization need not be monotone because particular future observations can alter the smoother posterior in different ways.

The information-theoretic identity \eqref{eq:expected_fingerprint_cmi} also establishes a direct connection to transfer entropy~\cite{schreiber2000measuring} and its generalization to causal networks, causation entropy~\cite{sun2014causation,sun2015causal}. Formally adapting the definition of causation entropy to the present setting, denoted by $\mathcal{C}$, we have
\begin{align}\label{eq:sec2_CE}
\mathcal{C}_{\mathbf{Y}_t\to\mathbf{X}_{(t,t+\tau]}\mid \mathbf{X}_{[0,t]}}&:=H\left(\mathbf{X}_{(t,t+\tau]}\mid\mathbf{X}_{[0,t]}\right)
-
H\left(\mathbf{X}_{(t,t+\tau]}\mid\mathbf{X}_{[0,t]},\mathbf{Y}_t\right)\notag\\
&=
I\left(\mathbf{Y}_t;\mathbf{X}_{(t,t+\tau]}\mid \mathbf{X}_{[0,t]}\right)
\equiv
\mathbb E_{\mathbf{X}_{[0,t+\tau]}}\bigl[\mathcal I(t,\tau)\bigr].
\end{align}
Therefore, the average causal fingerprint recovers the same conditional-information principle underlying causation entropy: a cause is informative when its knowledge reduces the uncertainty about a subsequent effect.
In addition, by symmetry of the conditional mutual information, the averaged fingerprint in \eqref{eq:sec2_CE} equals \(H(\mathbf Y_t\mid\mathbf X_{[0,t]})-H(\mathbf Y_t\mid\mathbf X_{[0,t+\tau]})\), expressing how subsequent observations reduce uncertainty about their earlier hidden cause. For an individual observed trajectory, the causal fingerprint provides a pathwise counterpart of this inverse-problem interpretation through Bayesian data assimilation.

For a fixed state time $t$, varying the future-window length $\tau$ shows how far the causal influence assigned to the hidden state $\mathbf{Y}_t$ remains visible in the subsequent observations $\mathbf{X}_{[0,t+\tau]}$. Moving the state time along a trajectory $\mathbf{X}$ shows how this influence changes as an event develops. The observed trajectory itself is held fixed, but the inferred hidden state and its fingerprint depend on the candidate model. Comparing fingerprints obtained from the same observations under different models therefore reveals how competing hidden mechanisms explain the same observed evolution.

\section{Mechanism Discrimination in Linear Gaussian Partially-Observed Systems}
\label{sec:linear}

The framework developed in Section~\ref{sec:theory} leaves open whether the temporal structure of the causal fingerprint can discriminate misrepresented hidden mechanisms or omitted latent components. In this section, we study this question in linear Gaussian systems, where the filtering and smoothing distributions, as well as the fingerprint itself, can be computed analytically in closed form. In this setting, we analyze the memory effects in the evolution of the observed variable and derive a characteristic timescale for the temporal reach of hidden causal influence.

We perform two tests in this section. In each test, a {perfect model} serves as the generator that produces a fixed reference trajectory of \(X\) under specific latent dynamical mechanisms. The perfect and surrogate models then serve as inference models for computing fingerprints along this same realized trajectory. The corresponding {surrogate models} omit or suppress mechanisms present in the perfect model but are nevertheless calibrated to reproduce the same selected observable statistics. We then investigate whether causal fingerprints can reveal the mechanistic differences between the perfect and surrogate models despite this statistical agreement in the observable state space.

Test I considers a perfect model with both feedback-induced state memory and hidden-noise memory. Its two surrogates respectively remove the feedback-induced state memory and strongly suppress the hidden-noise memory. Test II considers a perfect model with slow and fast hidden forcing components and a lower-dimensional surrogate with just a single hidden forcing component.

\subsection{Memory effects and observable statistics}\label{sec:sec3_linear_systems}
In this subsection, we analyze memory effects in the observed variable of a linear Gaussian model and derive closed-form expressions for its stationary observable statistics. These results will be used in the tests in Sections~\ref{sec:linear-test-i}--\ref{sec:linear-test-ii}.

We consider a two-dimensional linear Gaussian model
\begin{subequations}\label{eq:sec3_perfect_model}
\begin{align}
dX_t &= (-aX_t + Y_t)dt + \sigma_XdW_t^X,
\label{eq:sec3_perfect_model_X}\\
dY_t &= (cX_t - bY_t)dt + \sigma_YdW_t^Y,
\label{eq:sec3_perfect_model_Y}
\end{align}
\end{subequations}
where $X_t \in \mathbb{R}$ is the observed variable and $Y_t \in \mathbb{R}$ is the hidden one. The term $Y_t$ acts as a hidden forcing on $X_t$, while $cX_t$ represents feedback from $X_t$ to $Y_t$. The parameters $a>0$ and $b>0$ are the damping rates of $X_t$ and $Y_t$, respectively, and $W_t^X$ and $W_t^Y$ are independent standard Wiener processes with $\sigma_X>0$ and $\sigma_Y>0$ being the corresponding noise amplitudes. We assume $a>0$, $b>0$, and $ab-c>0$, so the system is stable and admits a stationary Gaussian distribution. We also assume that the initial state \((X_0,Y_0)\) follows a joint Gaussian distribution.

\paragraph{Memory effects decomposition.}
Under \eqref{eq:sec3_perfect_model_Y}, the hidden variable \(Y_t\) retains the effects of past observed states and past hidden stochastic forcing and carries them into the current evolution of \(X_t\) in \eqref{eq:sec3_perfect_model_X}. We refer to this dependence of the observed dynamics on past states and forcing as {memory effects}~\cite{zwanzig1961memory}. To make these memory effects explicit, we first solve the hidden-state equation \eqref{eq:sec3_perfect_model_Y}, obtaining \cite{sarkka2019applied}:
\begin{equation}
Y_t
=
e^{-bt}Y_0
+
c\int_0^t e^{-b(t-s)} X_sds
+
\sigma_Y \int_0^t e^{-b(t-s)}dW_s^Y.
\label{eq:sec3_Y_solution}
\end{equation}
Substituting this solution into \eqref{eq:sec3_perfect_model_X} to eliminate $Y_t$ yields
\begin{equation}
dX_t
=
\left[
-aX_t
+
c\int_0^t e^{-b(t-s)}X_sds
+
e^{-bt}Y_0
+
\sigma_Y \int_0^t e^{-b(t-s)}dW_s^Y
\right]dt
+ \sigma_XdW_t^X.
\label{eq:sec3_X_memory}
\end{equation}

Equation \eqref{eq:sec3_X_memory} decomposes these memory effects
into two distinct contributions stemming from the hidden state. The
term \(c\int_0^t e^{-b(t-s)}X_s\,ds\) is a {feedback-induced state memory}:
past values of \(X\) influence the hidden state \(Y_t\) through feedback forcing, and \(Y_t\) in turn drives the current evolution of \(X_t\). The stochastic
convolution \(\sigma_Y\int_0^t e^{-b(t-s)}\,dW_s^Y\) is a
{hidden-noise memory}, which carries the influence of past hidden
stochastic forcing into the current dynamics of \(X_t\). In addition, \(e^{-bt}Y_0\) represents the decaying influence of the initial hidden state, while \(\sigma_X\,dW_t^X\) is direct white-noise forcing.

\paragraph{Stationary observable statistics.}
We next consider the system \eqref{eq:sec3_perfect_model} in its stationary regime or equilibrium (i.e., for the following, we let $t\to\infty$) and summarize its observable statistics: the mean and variance of the stationary distribution of \(X_t\), and a decorrelation time that characterizes the temporal persistence of \(X\) \cite{gardiner2009stochastic, sarkka2019applied, majdaharlim2012filtering}. A derivation of the variance and decorrelation-time formulas is given in Appendix~\ref{app:stationary_statistics}.

Because the drift in \eqref{eq:sec3_perfect_model} contains no constant forcing term,
the stationary means are simply
\begin{equation}
\mathbb{E}[X_t]=\mathbb{E}[Y_t]=0,
\label{eq:sec3_calibration_mean}
\end{equation}
while the stationary variance of \(X_t\) is given by
\begin{equation}
V_X:=\operatorname{Var}(X_t)
=
\frac{\sigma_Y^2+\bigl(b(a+b)-c\bigr)\sigma_X^2}{2(a+b)(ab-c)}.
\label{eq:sec3_calibration_cov}
\end{equation}
Define now the {normalized observable autocorrelation} by
\[
\rho_X(\tau)
:=
\frac{\operatorname{Cov}(X_{t+\tau},X_t)}{V_X}.
\]
Here, \(\tau\) is the time lag; at stationarity, this correlation depends only on \(\tau\), and \(\rho_X(0)=1\). We define the observable {decorrelation time} using the integral of the squared autocorrelation~\cite{delsole2001persistent}:
\begin{equation}
\mathcal T_X
:=
\int_0^\infty \rho_X(\tau)^2\,d\tau
=
\frac{m_X^2+2bm_X+b(a+b)-c}{2(a+b)(ab-c)},
\label{eq:sec3_calibration_tcorr2}
\end{equation}
where \(m_X:=\operatorname{Cov}(Y_t,X_t)/V_X=a-\sigma_X^2/(2V_X)\) is the ratio of the stationary cross-covariance to the stationary observable variance.
Larger values of \(\mathcal T_X\) indicate that the autocorrelation retains a large magnitude over a longer range of time lags. Squaring \(\rho_X(\tau)\) in \eqref{eq:sec3_calibration_tcorr2} avoids cancellation between positive and negative autocorrelations in the integral for oscillatory dynamics.

\subsection{Analytical evaluation and characteristic timescale of the fingerprint}\label{sec:sec3_closed_form}
In this subsection, we give an explicit expression for the causal fingerprint in the linear Gaussian model and derive a characteristic timescale for the temporal reach of hidden causal influence.

For the scalar hidden state \(Y_t\) in \eqref{eq:sec3_perfect_model}, the filtering and smoothing distributions in \eqref{eq:info_gain_natural} are one-dimensional Gaussians~\cite{sarkka2013bayesian, liptser2001statistics,liptser2001applications}:
\[
p_{\mathrm f}(Y_t;t)
:=
p(Y_t\mid X_{[0,t]})
=
\mathcal N\bigl(\mu_{\mathrm f}(t),R_{\mathrm f}(t)\bigr),
\qquad
p_{\mathrm s}(Y_t;t,\tau)
:=
p(Y_t\mid X_{[0,t+\tau]})
=
\mathcal N\bigl(\mu_{\mathrm s}(t,\tau),R_{\mathrm s}(t,\tau)\bigr).
\]
Here, \(\mu_{\mathrm f}\) and \(\mu_{\mathrm s}\) are the filtering and smoothing means, respectively, and \(R_{\mathrm f}\) and \(R_{\mathrm s}\) are the corresponding variances. For these Gaussian distributions, the fingerprint definition \eqref{eq:info_gain_natural} takes the following closed form, known as the signal-dispersion formula \cite{cai2002mathematical}:
\begin{equation}
\mathcal I(t,\tau)
=
\frac{1}{2}
\left[
\frac{\bigl(\mu_{\mathrm s}(t,\tau)-\mu_{\mathrm f}(t)\bigr)^2}
{R_{\mathrm f}(t)}
+\frac{R_{\mathrm s}(t,\tau)}{R_{\mathrm f}(t)}
-1
-\log\frac{R_{\mathrm s}(t,\tau)}{R_{\mathrm f}(t)}
\right].
\label{eq:sec3_gaussian_kl}
\end{equation}
Once the filtering and smoothing means and variances have been computed, they can be substituted into \eqref{eq:sec3_gaussian_kl} to evaluate the fingerprint.

\paragraph{Filtering variance and online smoothing.}
To analyze how the smoothing distribution at a fixed state time \(t\) changes as future observations are added (i.e., by varying $\tau$), we use an {online smoother}~\cite{andreou2026adaptive}. The filtering variance and online smoothing equations below follow from the general formulation for a {conditional Gaussian nonlinear system} (CGNS) \cite{chenmajda2018conditional}, as outlined in Appendix~\ref{app:cgns_derivations}, which also provides the full filtering equations and the equivalent {fixed-interval smoothing} formulation~\cite{liptser2001statistics,liptser2001applications, chenmajda2016filtering, chenmajda2020smoothing, andreou2024martingale}.

The filtering variance satisfies a Riccati differential equation:
\begin{equation}
\mathrm{d} R_{\mathrm f}(t)
=
\left[
-2bR_{\mathrm f}(t)
+\sigma_Y^2
-\frac{R_{\mathrm f}(t)^2}{\sigma_X^2}
\right]\mathrm{d} t.
\label{eq:sec3_scalar_filter_variance}
\end{equation}
For any nonnegative initial variance \(R_{\mathrm f}(0)\), the solution of \eqref{eq:sec3_scalar_filter_variance} converges exponentially to the steady value
\begin{equation}
R_*
=
\sigma_X^2
\left(
-b+\sqrt{b^2+\frac{\sigma_Y^2}{\sigma_X^2}}
\right).
\label{eq:sec3_scalar_steady_filter_covariance}
\end{equation}

At a fixed state time \(t\), the online smoother updates the smoothing mean \(\mu_{\mathrm s}(t,\tau)\) and variance \(R_{\mathrm s}(t,\tau)\) as the future observation window expands. These updates can be represented through the scalar {smoothing influence kernel or propagator} \(\Phi(t,\tau)\), defined by
\begin{equation}
\frac{\partial}{\partial\tau}\Phi(t,\tau)
=
-\left(
b+\frac{R_{\mathrm f}(t+\tau)}{\sigma_X^2}
\right)\Phi(t,\tau),
\qquad
\Phi(t,0)=1.
\label{eq:sec3_scalar_influence_kernel_evolution}
\end{equation}
With \(t\) held fixed, we use \(\mathrm{d}_\tau\) to denote the It\^o differential as \(\tau\) increases; in particular, \(\mathrm{d}_\tau X_{t+\tau}\) denotes the observation increment at the advancing endpoint \(t+\tau\).
This kernel controls how strongly a new observation at \(t+\tau\) influences the Bayesian update about the earlier hidden state \(Y_t\). Starting from
\(\mu_{\mathrm s}(t,0)=\mu_{\mathrm f}(t)\) and
\(R_{\mathrm s}(t,0)=R_{\mathrm f}(t)\), the smoothing mean evolves in $\tau$ according to
\begin{equation}
\mathrm{d}_\tau\mu_{\mathrm s}(t,\tau)
=
\frac{R_{\mathrm f}(t)\Phi(t,\tau)}{\sigma_X^2}
\left[
\mathrm{d}_\tau X_{t+\tau}
-
\left(
\mu_{\mathrm f}(t+\tau)-aX_{t+\tau}
\right)\mathrm{d}\tau
\right].
\label{eq:sec3_scalar_online_smoother_mean}
\end{equation}
Here, the bracketed term is the difference between the advanced observed increment and its prediction from the filtering estimate based on \eqref{eq:sec3_perfect_model_X}, known as the {filtering innovation}. The corresponding smoothing variance satisfies
\begin{equation}
\mathrm{d}_\tau R_{\mathrm s}(t,\tau)
=
-\frac{R_{\mathrm f}(t)^2\Phi(t,\tau)^2}{\sigma_X^2}\mathrm{d}\tau.
\label{eq:sec3_scalar_online_smoother_covariance}
\end{equation}

\paragraph{Characteristic fingerprint timescale.}
For a fixed state time \(t\) after the filtering variance has sufficiently approached \(R_*\), we may approximate \(R_{\mathrm f}(t+\tau)\) by \(R_*\). The influence kernel in \eqref{eq:sec3_scalar_influence_kernel_evolution} then depends only on the future-window length \(\tau\), with its solution given by:
\begin{equation}
\Phi(t,\tau)
=
\Phi(\tau)
=
\exp\left[
-\left(
b+\frac{R_*}{\sigma_X^2}
\right)\tau
\right].
\label{eq:sec3_scalar_influence_kernel}
\end{equation}
The corresponding characteristic timescale is therefore
\begin{equation}
\tau_*
:=
\left(b+\frac{R_*}{\sigma_X^2}\right)^{-1}
=
\left(b^2+\frac{\sigma_Y^2}{\sigma_X^2}\right)^{-1/2}.
\label{eq:sec3_scalar_influence_rate}
\end{equation}
Thus, \(\Phi(\tau)=e^{-\tau/\tau_*}\). The kernel decay rate \(\tau_*^{-1}=b+R_*/\sigma_X^2\) combines the intrinsic hidden-state decay \(b\) with the steady-state {Kalman-Bucy gain} \(R_*/\sigma_X^2\).
In \eqref{eq:sec3_scalar_online_smoother_mean}, \(\Phi\) weights the observation innovation in the smoothing mean update, while \(\Phi^2\) controls the variance reduction in \eqref{eq:sec3_scalar_online_smoother_covariance}. As these weights decay with $\tau$, more distant observations have a diminishing effect on the posterior smoother inference of \(Y_t\), and the smoothing mean and variance converge as \(\tau\to\infty\). The filtering moments \(\mu_{\mathrm f}(t)\) and \(R_{\mathrm f}(t)\) do not depend on \(\tau\) at this fixed \(t\), so the fingerprint in \eqref{eq:sec3_gaussian_kl} necessarily converges as well. We therefore refer to \(\tau_*\) as a {characteristic fingerprint timescale}. It provides an effective information horizon over which future observations materially refine the inference of \(Y_t\); within the model-based causal interpretation adopted here, this horizon quantifies the temporal reach of the hidden state's influence on \(X\).

\subsection{Test I: feedback-induced state memory versus hidden-noise memory}\label{sec:linear-test-i}

In Test I we study whether causal fingerprints can discriminate between models that retain, omit, or suppress different hidden memory effects while reproducing the same selected observable summary statistics. The data-generating perfect model \eqref{eq:sec3_perfect_model} contains both feedback-induced state memory and hidden-noise memory. The no-feedback surrogate removes the former, whereas the feedback-memory-dominated surrogate strongly suppresses the latter; both surrogates' model parameters are calibrated to reproduce the perfect model's stationary mean of \(\mathbb{E}[X_t]=0\), variance \(V_X\) (see \eqref{eq:sec3_calibration_cov}), and decorrelation time \(\mathcal T_X\) (see \eqref{eq:sec3_calibration_tcorr2}).
We compute the fingerprints of all three models from the same observed trajectory generated by the perfect model, examining how the different memory mechanisms affect the patterns of hidden causal influence.

\subsubsection{Perfect model and surrogates}

Equation \eqref{eq:sec3_X_memory} expresses the two memory effects of the dynamics of \(X_t\) as
\[
\begin{aligned}
\text{feedback-induced state memory:}\quad
&c\int_0^t e^{-b(t-s)}X_s\,ds,\\
\text{hidden-noise memory:}\quad
&\sigma_Y\int_0^t e^{-b(t-s)}dW_s^Y.
\end{aligned}
\]
The common kernel \(e^{-b(t-s)}\) determines how contributions from earlier times decay, while \(c\) and \(\sigma_Y\) weight the state and stochastic-forcing contributions, respectively. All three models have the two-dimensional linear Gaussian form \eqref{eq:sec3_perfect_model}, but represent these memory effects differently:

\paragraph{Perfect model.}
The data-generating perfect model has \(c\ne0\) and \(\sigma_Y>0\), retaining both memory contributions.

\paragraph{No-feedback surrogate.}
The first surrogate sets \(c=0\) while retaining \(\sigma_Y>0\). The feedback-induced state memory term vanishes, leaving only hidden-noise memory.

\paragraph{Feedback-memory-dominated surrogate.}
The second surrogate retains \(c\ne0\) but uses a small nonzero hidden-noise amplitude \(\sigma_Y\). Its feedback-induced state memory therefore remains active,
while hidden-noise memory is strongly suppressed.

\subsubsection{Experimental setup and parameter choices}

In Test I, we consider two parameter regimes of the perfect model: oscillatory and non-oscillatory. The deterministic drift matrix \(\mathbf A\) of \eqref{eq:sec3_perfect_model} and its eigenvalues are
\begin{equation}
\mathbf A:=\begin{pmatrix}-a&1\\c&-b\end{pmatrix},
\qquad
\lambda_{\pm}
=
\frac{-a-b\pm\sqrt{\Delta}}{2},
\qquad
\Delta:=(a-b)^2+4c.
\label{eq:sec3_timescales_eigenvalues}
\end{equation}
We call the deterministic dynamics {oscillatory} when \(\Delta<0\), and {non-oscillatory} when \(\Delta>0\). The perfect-model parameters used in Test I and the corresponding stationary observable statistics evaluated from \eqref{eq:sec3_calibration_mean}--\eqref{eq:sec3_calibration_tcorr2} are listed in Table~\ref{tab:testI_perfect_parameters_statistics}.
\begin{table}[H]
\centering
\caption{Parameters and stationary observable statistics of the perfect model in Test I.}
\label{tab:testI_perfect_parameters_statistics}
\renewcommand{\arraystretch}{1.15}
\begin{ruledtabular}
\begin{tabular}{ccccccccc}
regime\hspace{0.5em} & $a$ & $b$ & $c$ & $\sigma_X$ & $\sigma_Y$\hspace{0.5em} & $\mathbb{E}[X_t]$ & $\mathrm{Var}(X_t)$ & $\mathcal T_X$\\[2pt]
\hline
oscillatory\hspace{0.5em} & $0.12$ & $0.15$ & $-3.0$ & $0.45$ & $0.10$\hspace{0.5em} & $0$ & $0.384$ & $1.852$\\
non-oscillatory\hspace{0.5em} & $1.30$ & $1.20$ & $0.34$ & $0.25$ & $0.08$\hspace{0.5em} & $0$ & $0.028$ & $0.519$\\
\end{tabular}
\end{ruledtabular}
\end{table}

For simplicity, we hold \(\sigma_X\) at the corresponding perfect-model value in both surrogates. In practice, \(\sigma_X\) can be estimated from the quadratic variation of the observed process \(X\)~\cite{sarkka2019applied}.
The stationary mean is zero for all three models. We select the remaining surrogate parameters to match the perfect model's stationary mean, variance \(V_X\), and decorrelation time \(\mathcal T_X\) exactly. Within each regime, we first choose a common value for the hidden-state damping rate \(b\) for the two surrogates. We then solve \eqref{eq:sec3_calibration_cov} and \eqref{eq:sec3_calibration_tcorr2} for the remaining degrees of freedom: \((a,\sigma_Y)\) with \(c=0\) in the no-feedback surrogate, or for \((a,c)\) with a small prescribed nonzero \(\sigma_Y\) in the feedback-memory-dominated surrogate. Thus, these surrogates are not single-parameter ablations: after one memory mechanism is removed or strongly suppressed, the remaining parameters are recalibrated to preserve the selected observable statistics. The test asks whether the fingerprint can still reveal the resulting mechanistic difference after this statistical compensation. The resulting surrogate model parameters used in Test I are listed in Table~\ref{tab:testI_surrogate_parameters}.
\begin{table}[H]
\centering
\caption{Parameters of the surrogate models in Test I.}
\label{tab:testI_surrogate_parameters}
\renewcommand{\arraystretch}{1.15}
\begin{ruledtabular}
\begin{tabular}{ccccccc}
perfect-model regime & surrogate\hspace{0.5em} & $a$ & $b$ & $c$ & $\sigma_X$ & $\sigma_Y$\\[2pt]
\hline
oscillatory & no-feedback\hspace{0.5em} & $0.271$ & $3$ & $0$ & $0.45$ & $0.238$\\
oscillatory & $\vcenter{\hbox{\shortstack{feedback-memory-\\dominated}}}$\hspace{0.5em} & $0.226$ & $3$ & $-0.125$ & $0.45$ & $0.014$\\[3pt]
\hline
non-oscillatory & no-feedback\hspace{0.5em} & $1.673$ & $0.6$ & $0$ & $0.25$ & $0.210$\\
non-oscillatory & $\vcenter{\hbox{\shortstack{feedback-memory-\\dominated}}}$\hspace{0.5em} & $1.209$ & $0.6$ & $0.160$ & $0.25$ & $0.008$\\[3pt]
\end{tabular}
\end{ruledtabular}
\end{table}

For each of the oscillatory and non-oscillatory regimes, we generate and fix one reference trajectory \(X_{[0,50]}\) from the corresponding perfect model. Using the same regime-specific trajectory, we evaluate the fingerprints of the perfect model and the two surrogates over \(0\le t\le25\) and \(0\le\tau\le25\). The evaluation uses the filtering moments from \eqref{eq:sec3_scalar_filter_mean} and \eqref{eq:sec3_scalar_filter_variance}, the fixed-interval smoothing moments from \eqref{eq:sec3_scalar_backward_smoother}, and the Gaussian KL expression \eqref{eq:sec3_gaussian_kl}. Unless otherwise stated, all numerical simulations in this paper use the Euler-Maruyama scheme~\cite{sarkka2019applied} with time step \(\Delta t=10^{-3}\).

\subsubsection{Results}

\begin{figure}[tbp!]
\centering
\includegraphics[width=1.0\textwidth]{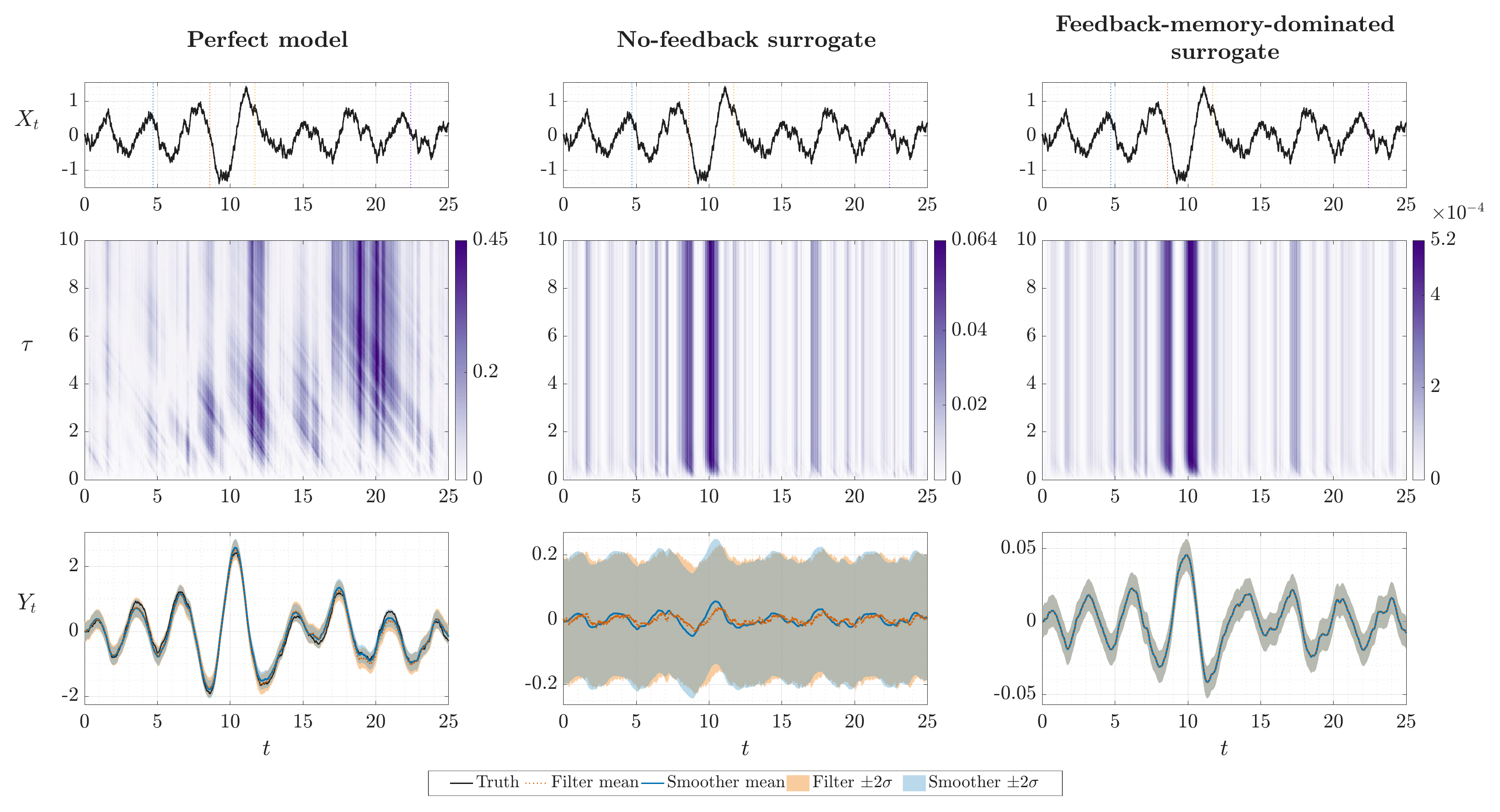}
\caption{Reference trajectory, causal fingerprints, and hidden-state estimates for the oscillatory regime in Test I. Columns show the perfect model, the no-feedback surrogate, and the feedback-memory-dominated surrogate. The top row repeats the same perfect-model reference trajectory \(X\) in all three columns. The middle row shows \(\mathcal I(t,\tau)\) over \(0\le t\le25\) and \(0\le\tau\le10\); the bottom row shows the posterior estimates of \(Y_t\). Orange dotted and blue solid curves show the filtering and 25-time-unit fixed-lag smoothing means, with shaded bands indicating \(\pm2\) posterior standard deviations. The solid black curve in the perfect-model panel shows the true hidden trajectory. Throughout this paper, each fingerprint heat map uses its own color scale, clipped at the 99.5th percentile. Vertical dotted lines in the top row mark the state times selected for the fingerprint curves in Figure~\ref{fig:sec3_testI_curves}, with corresponding colors.}
\label{fig:sec3_testI_oscillatory_combined}
\end{figure}

\begin{figure}[tbp!]
\centering
\includegraphics[width=1.0\textwidth]{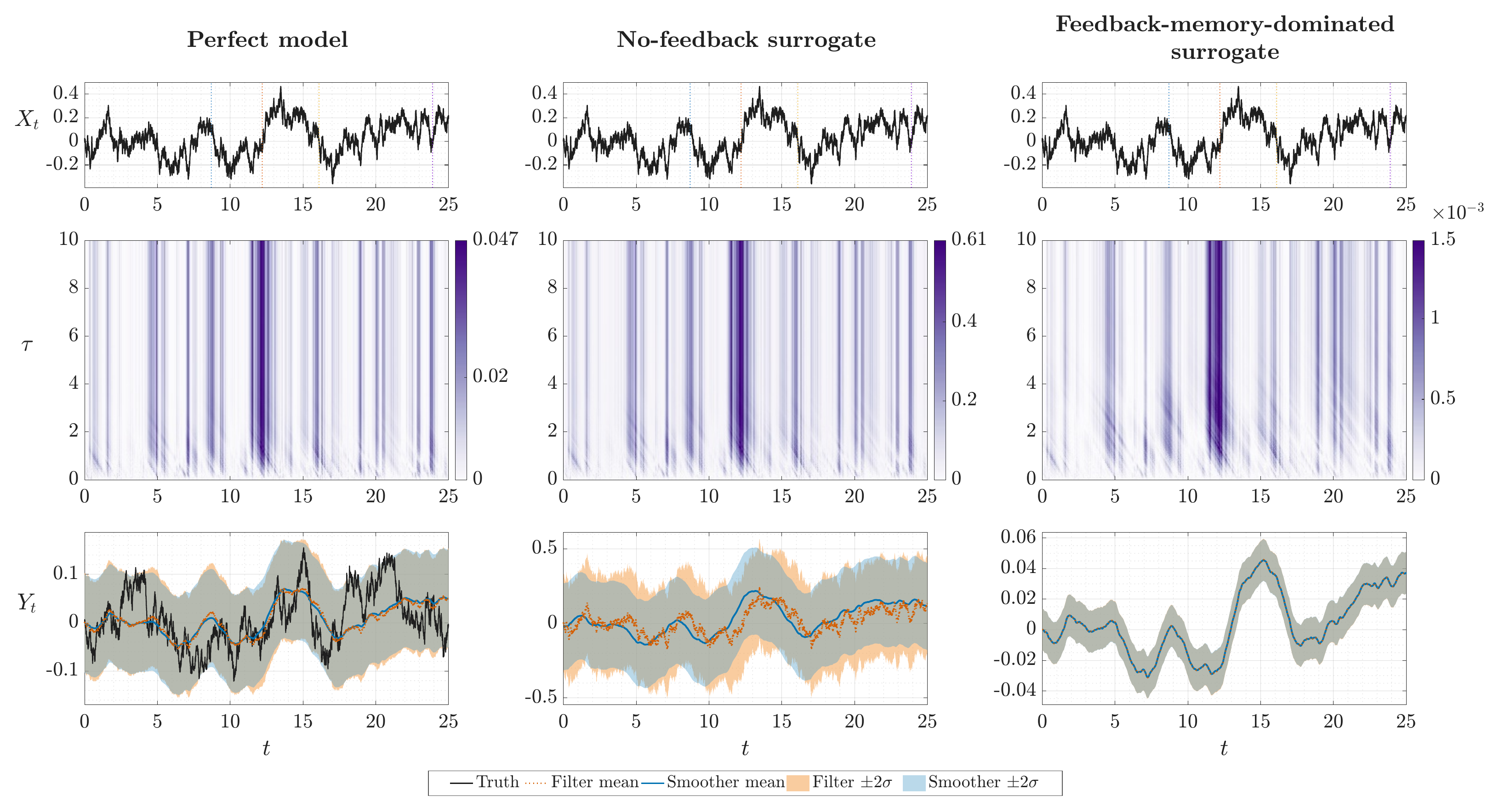}
\caption{Similar to Figure~\ref{fig:sec3_testI_oscillatory_combined} but for the non-oscillatory regime in Test I.}
\label{fig:sec3_testI_nonoscillatory_combined}
\end{figure}

\begin{figure}[tbp!]
\centering
\includegraphics[width=1.0\textwidth]{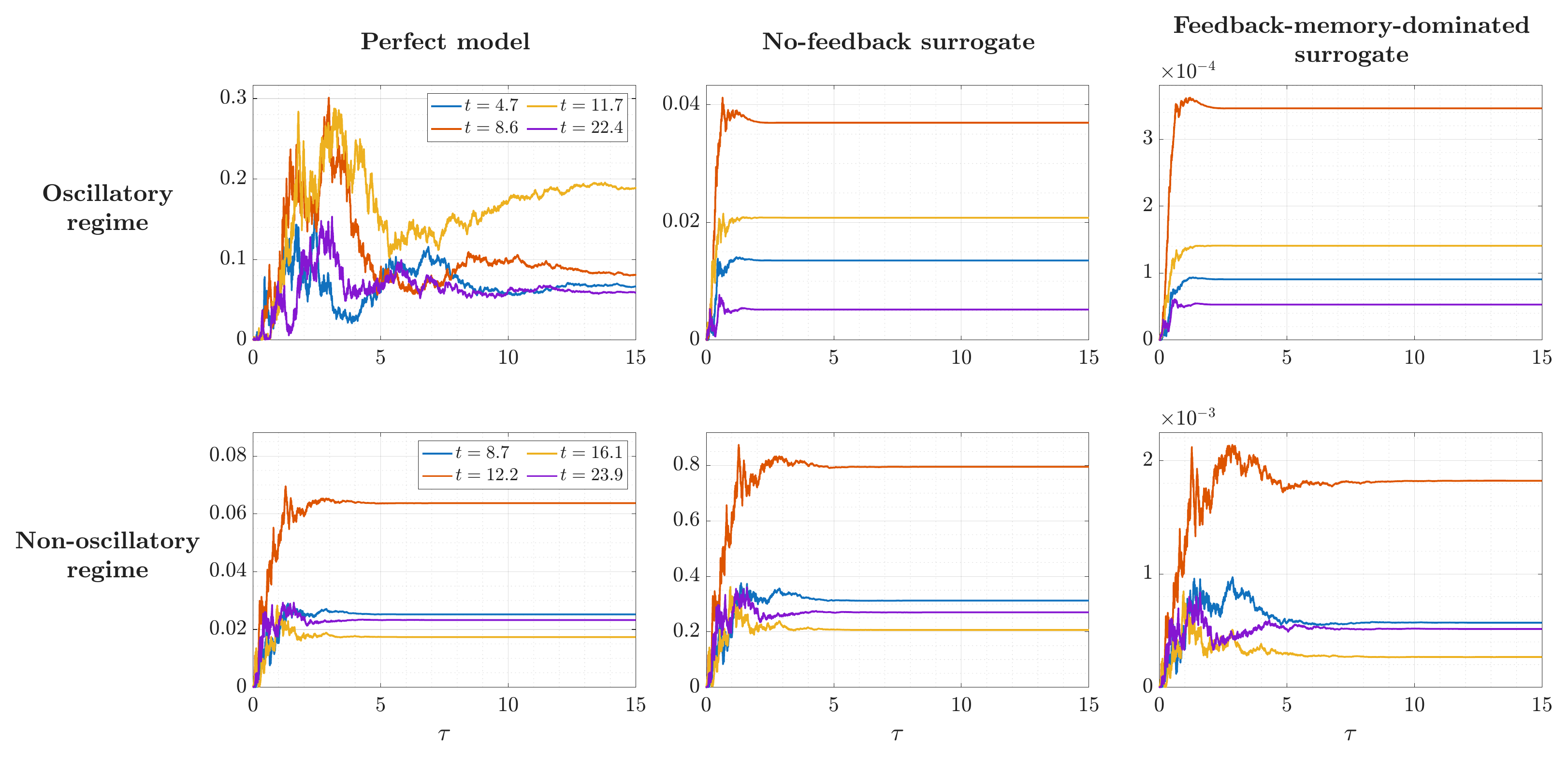}
\caption{Fixed-time fingerprint curves \(\mathcal I(t,\tau)\) for Test I, plotted against future-window length \(0\le\tau\le15\). Rows show the oscillatory (top) and non-oscillatory (bottom) perfect-model regimes. Columns show the perfect model, the no-feedback surrogate, and the feedback-memory-dominated surrogate. Curves are evaluated at state times \(t=4.7,8.6,11.7,22.4\) in the top row and \(t=8.7,12.2,16.1,23.9\) in the bottom row, with consistent colors across all three models within each row. These state times are marked by vertical dotted lines of corresponding colors in the top rows of Figures~\ref{fig:sec3_testI_oscillatory_combined} and~\ref{fig:sec3_testI_nonoscillatory_combined}, respectively.}
\label{fig:sec3_testI_curves}
\end{figure}

Figures~\ref{fig:sec3_testI_oscillatory_combined} and~\ref{fig:sec3_testI_nonoscillatory_combined} show the reference observed trajectory of \(X\) and the corresponding causal fingerprint heat maps for the perfect and surrogate models, across state times \(t\) and future-window lengths \(\tau\), for the oscillatory and non-oscillatory regimes, respectively. In the oscillatory regime, the lower part of the perfect-model fingerprint heat map shows prominent diagonal streaks extending upward and to the left (column 1, row 2, Figure~\ref{fig:sec3_testI_oscillatory_combined}). These streaks, represented by \(t+\tau=s^\star\) lines, have \(t\)-axis intercepts that correspond to peaks, troughs, and sharp turns in \(X\), for example at \(s^\star\approx9.2\), \(11\), and \(17\). Each streak traces how hidden states at successively earlier times contribute to the same observed event at \(s^\star\). Figure~\ref{fig:sec3_testI_curves} shows fixed-time fingerprint curves as functions of \(\tau\) at selected fixed state times \(t\). For each model, the curves at different fixed state times \(t\) approach their corresponding plateaus over comparable future-window lengths. Consistent with the analysis in Section~\ref{sec:sec3_closed_form}, this occurs after several multiples of that model's characteristic fingerprint timescale \(\tau_*\). This timescale therefore marks an effective information horizon for the influence of \(Y_t\) on future \(X\): once the curves approach their plateaus, more distant observations provide little additional information for inferring \(Y_t\). In the oscillatory regime, the perfect-model curves exhibit oscillatory variability before approaching a plateau, reflecting the model's intrinsic oscillatory dynamics. In the non-oscillatory regime, in contrast, the perfect-model heat map lacks the prominent diagonal streaks seen in the oscillatory regime, and its fixed-time curves approach their plateaus over shorter future windows. For these parameter choices, hidden influence persists over longer future windows in the oscillatory regime. In the non-oscillatory regime, the strongest responses are instead associated more locally with individual noise-driven fluctuations.

Although the surrogates reproduce the same selected observable summary statistics as the perfect model, their causal fingerprints differ from the perfect-model fingerprint in structure. In the oscillatory regime, the characteristic fingerprint timescales of the no-feedback and feedback-memory-dominated surrogates are \(\tau_*=0.328\) and \(0.333\), respectively, much shorter than the perfect-model's value of \(3.73\) and its oscillation period \(4\pi/\sqrt{-\Delta}\approx3.63\). Consistently, their fixed-time curves approach a plateau before capturing the observed oscillatory pattern, and their heat maps do not reproduce the diagonal  streaks. The locations of the prominent fingerprint responses along the state-time axis also differ between models. In the middle row of Figure~\ref{fig:sec3_testI_oscillatory_combined}, these responses appear near \(t=11.7\) and \(t\approx20\) for the perfect model, but near \(t\approx8.6\) and \(10\) for the two surrogates, indicating that omitting or suppressing memory effects changes when the models infer strong hidden causal influence along the same observed trajectory. In the non-oscillatory regime, in contrast, the fingerprint heat maps in Figure~\ref{fig:sec3_testI_nonoscillatory_combined} show broadly similar patterns across the three models, with prominent fingerprint responses occurring at similar state times \(t\) despite differences in amplitude.

The two surrogates also produce markedly different fingerprint amplitudes. In particular, the feedback-memory-dominated surrogate produces a much smaller fingerprint than the no-feedback surrogate in both regimes. In this surrogate, weak hidden noise gives the hidden dynamics a high signal-to-noise ratio: \(Y_t\) is largely determined by past observations through feedback, so future observations add little information and the fingerprint remains small. Consistently, its filtering and smoothing means and uncertainty bands nearly coincide in the bottom-right panel of Figure~\ref{fig:sec3_testI_nonoscillatory_combined}. Together, these results show that causal fingerprints reveal changes in the strength and temporal structure of hidden causal influence when memory mechanisms are omitted or suppressed.

\subsection{Test II: one-timescale versus two-timescale hidden forcing}\label{sec:linear-test-ii}
In Test II we investigate whether causal fingerprints can discriminate the omission of a hidden forcing timescale in a lower-dimensional reduced-order model. The perfect model is temporally multiscale and drives \(X_t\) with slow and fast hidden forcing components, whereas the reduced model uses a single hidden forcing component. Similar to Test I, both the perfect and reduced models are calibrated to produce the same selected stationary statistics (the mean, variance \(V_X\), and decorrelation time \(\mathcal T_X\)) of \(X\), while their fingerprints are computed from the same reference trajectory of \(X\), generated by the perfect model.

Under the perfect model, the {joint causal fingerprint} of the slow-fast latent scales characterizes the joint causal influence of the two hidden forcing components on \(X\), while the {marginal causal fingerprints} characterize the causal influence of each component on \(X\). Test II examines whether the perfect-model fingerprint reveals distinct temporal patterns associated with the slow and fast hidden forcing components and whether the reduced-model fingerprint preserves this multiscale structure.

\subsubsection{Models and experimental setup}

\paragraph{Perfect model.}
The perfect model consists of one observed variable \(X_t\) driven by two hidden forcing components \(Y_t^{(1)}\) and \(Y_t^{(2)}\) with different damping timescales:
\begin{equation}
\begin{aligned}
dX_t &= (-aX_t + Y_t^{(1)} + Y_t^{(2)})dt + \sigma_XdW_t^X, \\
dY_t^{(1)} &= -b_1 Y_t^{(1)}dt + \sigma_Y^{(1)}dW_t^{(1)}, \\
dY_t^{(2)} &= -b_2 Y_t^{(2)}dt + \sigma_Y^{(2)}dW_t^{(2)}, \qquad b_1\ll b_2.
\end{aligned}
\label{eq:sec3_multiscale_perfect}
\end{equation}
Here, \(a>0\) is the damping rate of \(X_t\), \(b_1,b_2>0\) are the damping rates of the two hidden components, and \(\sigma_X,\sigma_Y^{(1)},\sigma_Y^{(2)}>0\) are the respective noise amplitudes, with \(W_t^X\), \(W_t^{(1)}\), and \(W_t^{(2)}\) being independent standard Wiener processes. We set \(b_1\ll b_2\), so \(Y_t^{(1)}\) is a slow variable with damping time \(b_1^{-1}\) and \(Y_t^{(2)}\) is a fast variable with damping time \(b_2^{-1}\); for sufficiently large $b_2$, $Y_t^{(2)}$ acts as a rapidly decorrelating forcing, while $Y_t^{(1)}$ retains colored-noise memory and influences $X_t$ over much longer times.

\paragraph{Reduced model.}
The reduced-order model replaces $(Y_t^{(1)},Y_t^{(2)})$ with a single hidden forcing \(Y_t\):
\begin{equation}
\begin{aligned}
dX_t &= (-aX_t + Y_t)dt + \sigma_XdW_t^X, \\
dY_t &= -bY_tdt + \sigma_YdW_t^Y.
\end{aligned}
\label{eq:sec3_multiscale_reduced}
\end{equation}
Here, \(a,b,\sigma_X,\sigma_Y>0\), and \(W_t^X\) and \(W_t^Y\) are independent standard Wiener processes. For our study, we set \(b=b_2\), so the reduced model essentially retains only the fast hidden forcing timescale and omits the slow one.

\paragraph{Parameter choices.}
We use the same observed-noise amplitude \(\sigma_X=0.25\) in both models. The reduced model's remaining free parameters \((a,\sigma_Y)\) are chosen to reproduce the perfect model's stationary variance and decorrelation time.
The parameters used in Test II are:
\[
\begin{array}{@{}l@{\quad}c@{\quad}l@{}}
\text{perfect model}
&:
&a=0.45,\;b_1=0.05,\;b_2=3.0,\;\sigma_X=0.25,\;\sigma_Y^{(1)}=0.034,\;\sigma_Y^{(2)}=0.683,
\\[0.5em]
\text{reduced model}
&:
&a=0.193,\;b=3.0,\;\sigma_X=0.25,\;\sigma_Y=0.169.
\end{array}
\]
With these parameters, both models have zero stationary mean and the same stationary variance and decorrelation time:
\[
\mathbb E[X_t]=0,
\qquad
\operatorname{Var}(X_t)=0.169,
\qquad
\mathcal T_X=2.602.
\]
The reduced-model statistics are computed from the analytic expressions in Section~\ref{sec:sec3_linear_systems}. The corresponding stationary variance and decorrelation-time formulas for the perfect model with two hidden forcing timescales are given in Appendix~\ref{app:stationary_statistics}.

We generate one reference trajectory \(X_{[0,50]}\) from the perfect model and use it to compute the causal fingerprints under both models over \(0\le t\le25\) and \(0\le\tau\le25\). Similar to Test I, both models have the CGNS structure, and their filtering and smoothing distributions, as well as their fingerprints, are computed using the equations in Section~\ref{sec:sec3_closed_form} and Appendix~\ref{app:cgns_derivations}.

\subsubsection{Joint and marginal fingerprints}

To characterize the causal influence of the slow and fast hidden forcing components on the observable $X_t$, we define their joint and marginal causal fingerprints.

For the perfect model, we define the joint filtering and smoothing distributions
and their componentwise marginals, for \(i=1,2\), as
\begin{align}
p_{\mathrm f}(Y_t^{(1)},Y_t^{(2)};t)
&=p(Y_t^{(1)},Y_t^{(2)}\mid X_{[0,t]}),
&\quad
p_{\mathrm s}(Y_t^{(1)},Y_t^{(2)};t,\tau)
&=p(Y_t^{(1)},Y_t^{(2)}\mid X_{[0,t+\tau]}),
\label{eq:sec3_test2_joint_posteriors}\\
p_{\mathrm f}^{(i)}(Y_t^{(i)};t)
&=p(Y_t^{(i)}\mid X_{[0,t]}),
&\quad
p_{\mathrm s}^{(i)}(Y_t^{(i)};t,\tau)
&=p(Y_t^{(i)}\mid X_{[0,t+\tau]}).
\label{eq:sec3_test2_marginal_posteriors}
\end{align}
Integrating the joint filtering density \(p_{\mathrm f}(Y_t^{(1)},Y_t^{(2)};t)\) over \(Y_t^{(2)}\) gives the marginal \(p_{\mathrm f}^{(1)}(Y_t^{(1)};t)\), while integrating over \(Y_t^{(1)}\) gives \(p_{\mathrm f}^{(2)}(Y_t^{(2)};t)\). The smoothing marginals \(p_{\mathrm s}^{(i)}(Y_t^{(i)};t,\tau)\) are obtained in the same way from \(p_{\mathrm s}(Y_t^{(1)},Y_t^{(2)};t,\tau)\).
The perfect model \eqref{eq:sec3_multiscale_perfect} has the CGNS structure, so the joint filtering and smoothing posteriors in \eqref{eq:sec3_test2_joint_posteriors} are Gaussian. The equations in Appendix~\ref{app:cgns_derivations} give their mean vectors \(\boldsymbol{\mu}_{\mathrm f}(t)\) and \(\boldsymbol{\mu}_{\mathrm s}(t,\tau)\), and covariance matrices \(\mathbf R_{\mathrm f}(t)\) and \(\mathbf R_{\mathrm s}(t,\tau)\), respectively. The marginal posteriors in \eqref{eq:sec3_test2_marginal_posteriors} are therefore also Gaussian: their means are the corresponding components of these mean vectors, and their variances are the corresponding diagonal entries of these covariance matrices. Specifically, component \(i\) has filtering variance \([\mathbf R_{\mathrm f}(t)]_{ii}\) and smoothing variance \([\mathbf R_{\mathrm s}(t,\tau)]_{ii}\), for \(i=1,2\). The corresponding joint and marginal fingerprints are defined as
\begin{equation}
\mathcal I(t,\tau)
=D_{\mathrm{KL}}\!\left(p_{\mathrm s}\Vert p_{\mathrm f}\right),
\qquad
\mathcal I^{(i)}(t,\tau)
=D_{\mathrm{KL}}\!\left(p_{\mathrm s}^{(i)}\Vert p_{\mathrm f}^{(i)}\right),
\quad i=1,2.
\label{eq:sec3_test2_fingerprint_definitions}
\end{equation}
The joint causal fingerprint characterizes the joint causal influence of the two hidden
forcing components on future observations, whereas the marginal causal fingerprints
characterize the slow- and fast-component influences separately.

If the two hidden components are conditionally independent under both the filtering and smoothing distributions, the posterior densities factorize as
\begin{equation}
p_{\mathrm f}(Y_t^{(1)},Y_t^{(2)};t)
=p_{\mathrm f}^{(1)}(Y_t^{(1)};t)p_{\mathrm f}^{(2)}(Y_t^{(2)};t),
\qquad
p_{\mathrm s}(Y_t^{(1)},Y_t^{(2)};t,\tau)
=p_{\mathrm s}^{(1)}(Y_t^{(1)};t,\tau)p_{\mathrm s}^{(2)}(Y_t^{(2)};t,\tau).
\label{eq:sec3_test2_posterior_factorization}
\end{equation}
Substituting the factorizations in
\eqref{eq:sec3_test2_posterior_factorization} into the KL definitions
\eqref{eq:sec3_test2_fingerprint_definitions} gives the additive decomposition into the marginal fingerprints~\cite{cover2006elements}:
\begin{equation}
\mathcal I(t,\tau)
=
\mathcal I^{(1)}(t,\tau)+\mathcal I^{(2)}(t,\tau).
\label{eq:sec3_test2_fingerprint_additivity}
\end{equation}
To assess the accuracy of this decomposition in Test II, under the inference model in \eqref{eq:sec3_multiscale_perfect}, we define the residual fingerprint:
\begin{equation}
\Delta\mathcal I(t,\tau)
:=
\mathcal I(t,\tau)
-\mathcal I^{(1)}(t,\tau)
-\mathcal I^{(2)}(t,\tau).
\label{eq:sec3_test2_additivity_residual}
\end{equation}
When this residual is small relative to the joint fingerprint, the sum of the slow and fast marginal fingerprints provides a good approximation to the joint fingerprint.

\subsubsection{Results}

\begin{figure}[tbp!]
\centering
\includegraphics[width=1.0\textwidth]{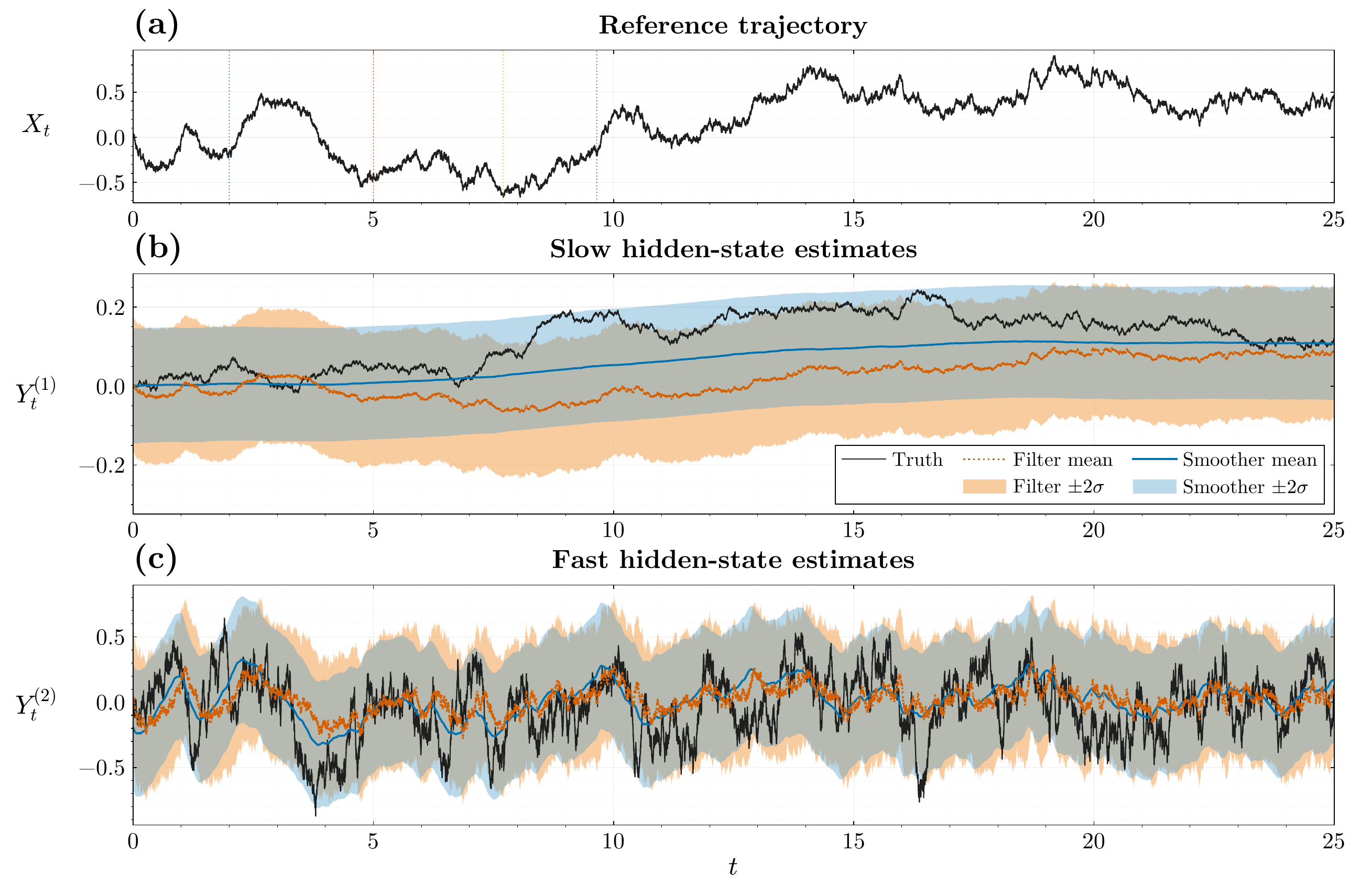}
\caption{Reference trajectory and hidden-state estimates under the perfect model in Test II. Panel (a) shows the reference trajectory \(X\) used to compute fingerprints under both the perfect and reduced models. Panels (b) and (c) show the slow and fast hidden components, \(Y_t^{(1)}\) and \(Y_t^{(2)}\), respectively. Solid black curves show the true hidden trajectories. Orange dotted and blue solid curves show the filtering and 25-time-unit fixed-lag smoothing means, with shaded bands indicating \(\pm2\) posterior standard deviations. Vertical dotted lines in panel (a) mark the state times selected for the fingerprint curves in Figure~\ref{fig:sec3_testII_perfect_decomposition}, with corresponding colors.}
\label{fig:sec3_testII_perfect_states}
\end{figure}

\begin{figure}[tbp!]
\centering
\includegraphics[width=1.0\textwidth]{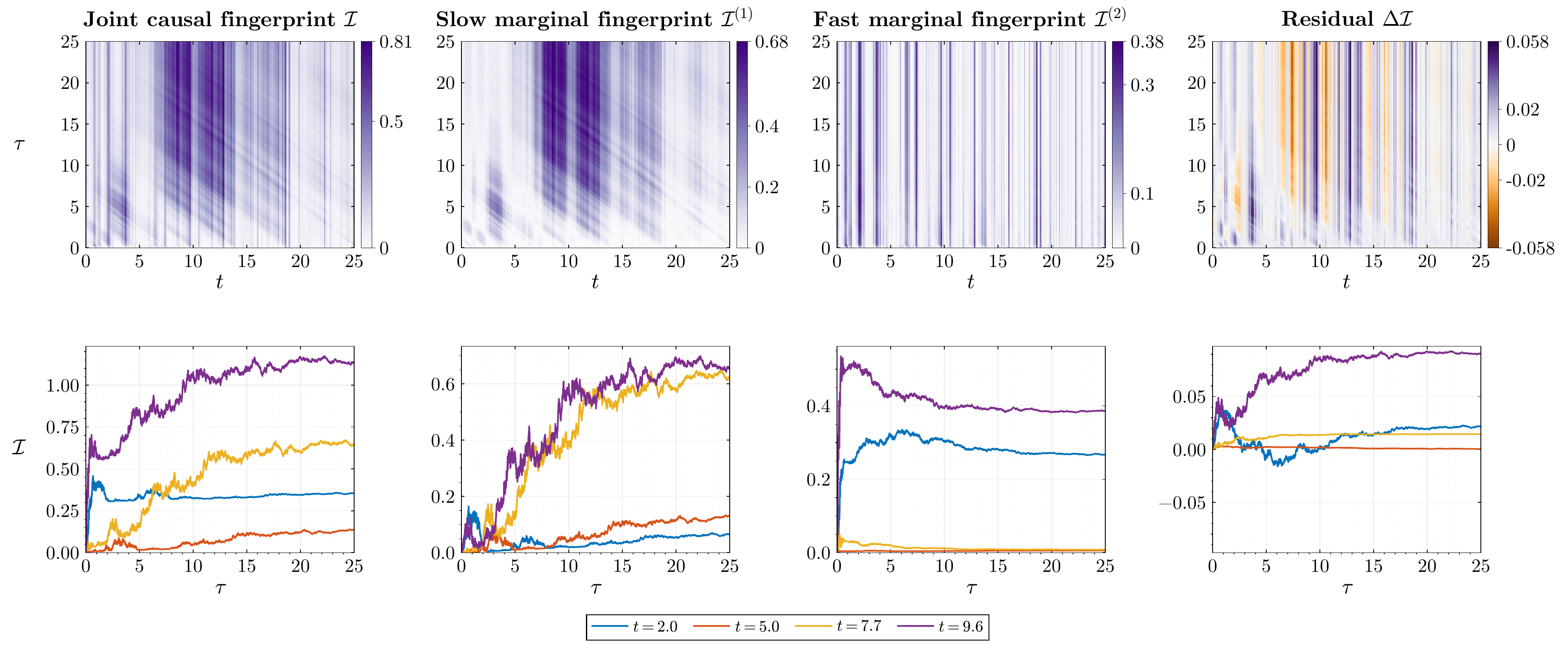}
\caption{Two-scale fingerprint decomposition for the perfect model in Test II. From left to right, columns show the joint fingerprint \(\mathcal I\), the slow and fast marginal fingerprints \(\mathcal I^{(1)}\) and \(\mathcal I^{(2)}\), and the signed residual \(\Delta\mathcal I=\mathcal I-\mathcal I^{(1)}-\mathcal I^{(2)}\). The top row shows heat maps over \(0\le t,\tau\le25\); the bottom row shows fixed-time fingerprint curves against \(\tau\). Curve colors identify the same state times \(t=2.0,5.0,7.7,9.6\) across all columns. These state times are marked by vertical dotted lines of corresponding colors in panel (a) of Figure~\ref{fig:sec3_testII_perfect_states}.}
\label{fig:sec3_testII_perfect_decomposition}
\end{figure}

Figure~\ref{fig:sec3_testII_perfect_states} shows the perfect-model reference
trajectory \(X\), the corresponding slow and fast hidden forcing paths
\(Y^{(1)}\) and \(Y^{(2)}\), and their filtering and smoothing posterior estimates.
Figure~\ref{fig:sec3_testII_perfect_decomposition} compares the joint
fingerprint, the two marginal fingerprints, and their residual. The residual
\(\Delta\mathcal I\) is small relative to the joint fingerprint, and has a
root-mean-square value of \(0.015\). Thus, \(\mathcal I^{(1)}+\mathcal I^{(2)}\) provides a good approximation to the joint fingerprint \(\mathcal I\). Consistent with this approximation, the maximum absolute posterior correlation between the two hidden components is \(0.084\) for filtering and \(0.151\) for smoothing, supporting approximate independence under both Gaussian posteriors. Here, the posterior correlation is computed as \(R_{12}/\sqrt{R_{11}R_{22}}\), with \(\mathbf R=\mathbf R_{\mathrm f}(t)\) for filtering and \(\mathbf R=\mathbf R_{\mathrm s}(t,\tau)\) for smoothing. In addition, the diagonal streaks in the joint fingerprint appear at corresponding locations in the slow marginal fingerprint, but are much less pronounced in the fast marginal fingerprint, indicating that the slow component accounts for much of this temporal structure.

The fixed-time fingerprint curves in the bottom row of Figure~\ref{fig:sec3_testII_perfect_decomposition} show distinct temporal patterns for the slow and fast marginal fingerprints. At \(t=2.0\) and \(t=9.6\), the fast marginal \(\mathcal I^{(2)}\) rises sharply within a short future window and then gradually approaches its plateau. The slow marginal \(\mathcal I^{(1)}\), most prominent at \(t=7.7\) and \(t=9.6\), develops gradually as the future observation window expands, indicating a longer temporal reach of the slow hidden component's influence on \(X\). Across the selected state
times, the joint fingerprint is dominated by the fast component at
\(t=2.0\) and the slow component at \(t=7.7\); both marginal fingerprints
are comparatively weak at \(t=5.0\) and strong at \(t=9.6\). The decomposition therefore reveals both the different temporal reaches of the two hidden components and their changing contributions at different state times \(t\).

\begin{figure}[tbp!]
\centering
\includegraphics[width=1.0\textwidth]{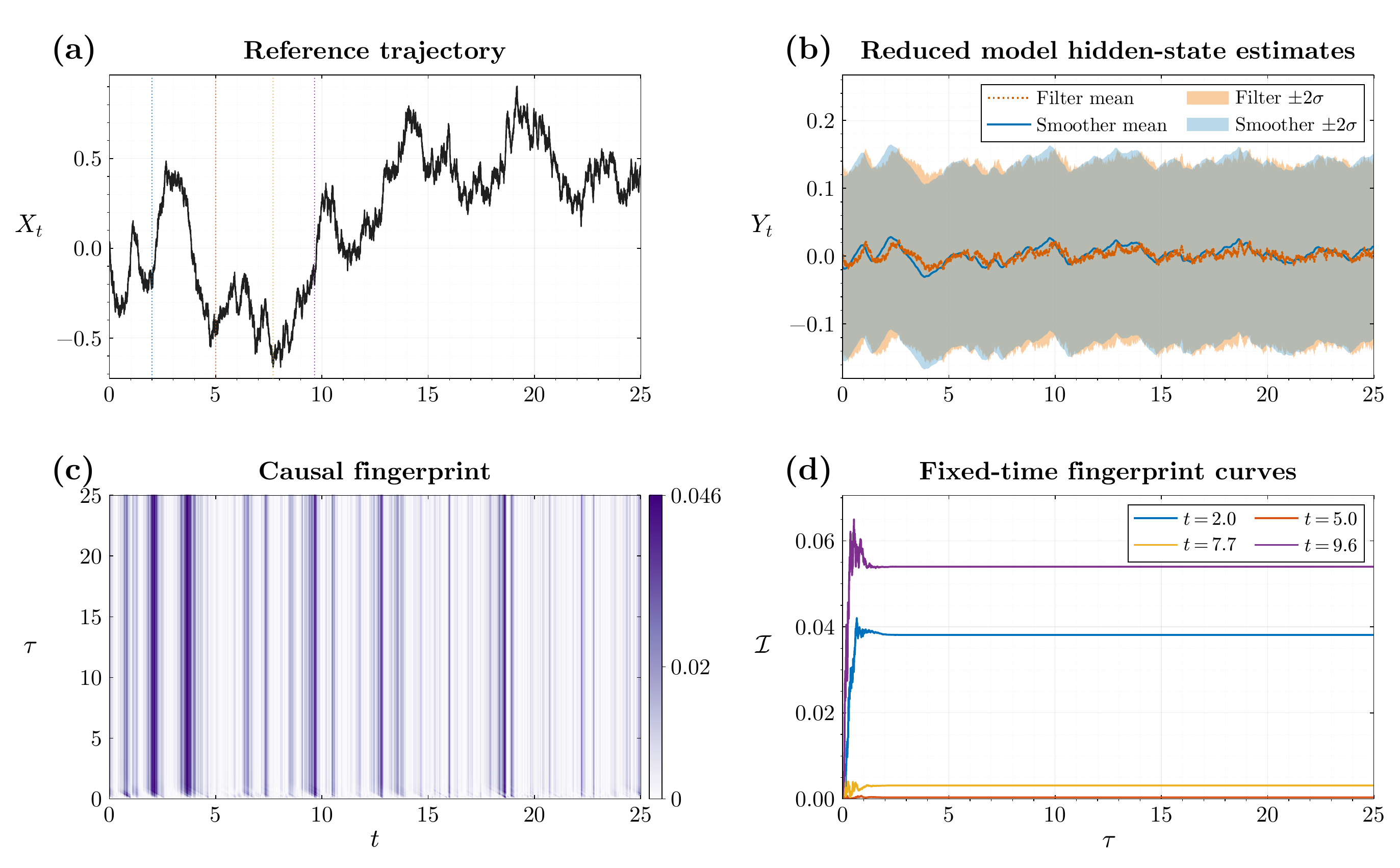}
\caption{Hidden-state estimates and causal fingerprint under the reduced model in Test II. Panel (a) shows the reference trajectory \(X\) generated by the perfect model. Panel (b) shows estimates of the reduced model's hidden state \(Y_t\): orange dotted and blue solid curves show the filtering and 25-time-unit fixed-lag smoothing means, with shaded bands indicating \(\pm2\) posterior standard deviations. Panel (c) shows \(\mathcal I(t,\tau)\) over \(0\le t,\tau\le25\). Panel (d) shows fixed-time fingerprint curves using the same four state times and colors as Figure~\ref{fig:sec3_testII_perfect_decomposition}. Vertical dotted lines in panel (a) mark the state times selected for the fingerprint curves in panel (d), with corresponding colors.}
\label{fig:sec3_testII_reduced_summary}
\end{figure}

Although the reduced model reproduces the same selected observable summary statistics as the perfect model, its scalar fingerprint in panels (c) and (d) of
Figure~\ref{fig:sec3_testII_reduced_summary} resembles the perfect model's
fast marginal fingerprint \(\mathcal I^{(2)}\), since we let $b=b_2$ in \eqref{eq:sec3_multiscale_reduced}, but lacks the response associated with the slow hidden forcing. At \(t=2.0\) and \(t=9.6\), the reduced-model fingerprint rises sharply over a short future window, similar to the perfect model's fast marginal fingerprint. However, it does not reproduce the gradual increase over longer future windows seen in the perfect model's slow marginal fingerprint at \(t=7.7\) and \(t=9.6\).

\section{Event-Level Causal Diagnostics for Nonlinear Mechanisms}
\label{sec:nonlinear}

In Section~\ref{sec:linear}, we studied how causal fingerprints respond when reduced models omit or suppress hidden mechanisms even after the surrogate models are calibrated to match the same selected observable summary statistics. We now extend this diagnostic approach to nonlinear examples, examining how competing model mechanisms assign hidden influence during different stages of the same observed event.

Models with different latent dynamical mechanisms can produce closely matched conventional summary statistics and similar classes of events, yet provide different mechanistic explanations for the same observed trajectory~\cite{aguilar2026limits, moser2026mechanisms}. In each case study, we consider two models that represent different mechanistic interpretations for the same type of observed event. The model parameters are chosen so that the two models produce closely matched values of selected observable statistics. With these parameters fixed, we generate one reference trajectory from each model containing the event of interest. Then, for each reference trajectory, both models are used as inference models: under this setup, the case study is {self-consistent} when the inference model is the generating model for the reference trajectory and {mechanism-mismatched} when it is the other model.

We consider three nonlinear examples: multiplicative instability versus additive hidden forcing as alternative mechanisms for intermittent extreme events; deterministic nonlinear feedback versus state-dependent stochastic forcing for oscillatory variability; and deterministic bistability versus hidden double-well forcing for regime switching and tipping.

\subsection{Extreme events: multiplicative instability versus additive forcing} \label{sec:extreme-events}

Intermittent extreme events can arise when a system temporarily enters an unstable regime or experiences strong additive forcing~\cite{mohamad2015probabilistic, chang2026extreme}. In this first example, we compare two models that produce bursts in the observable variable through these different mechanisms.

The first model represents {multiplicative instability}:
\begin{equation}
\begin{aligned}
dX_t &= (-aX_t+X_tY_t)dt + \sigma_XdW_t^X, \\
dY_t &= -bY_t\,dt + \sigma_YdW_t^Y .
\end{aligned}
\label{eq:sec4_mult_inst}
\end{equation}
In this nonlinear model, the hidden variable \(Y_t\) modulates the stability of \(X_t\) through the effective growth rate \(-a+Y_t\). When \(Y_t>a\), the observed dynamics enter an unstable regime. If this instability persists, exponential amplification can build up a large-amplitude burst in \(X_t\).

The competing model instead uses {additive hidden forcing}:
\begin{equation}
\begin{aligned}
dX_t &= (-aX_t+Y_t)dt + \sigma_XdW_t^X, \\
dY_t &= -bY_t\,dt + \sigma_YdW_t^Y .
\end{aligned}
\label{eq:sec4_additive}
\end{equation}
Both models use a hidden {Ornstein-Uhlenbeck} (OU) process \(Y_t\)~\cite{sarkka2019applied}; in the additive-forcing model, \(Y_t\) directly drives local changes in \(X_t\). A sufficiently strong and persistent forcing can produce large-amplitude bursts.

The parameters for the multiplicative-instability and additive-forcing models and the corresponding empirical estimates of the stationary mean, variance, and decorrelation time of \(X\) are listed in Table~\ref{tab:extreme_parameters}.
\begin{table}[H]
\centering
\caption{Parameters and empirical stationary observable statistics of the multiplicative-instability and additive-forcing models.}
\label{tab:extreme_parameters}
\renewcommand{\arraystretch}{1.15}
\begin{ruledtabular}
\begin{tabular}{cccccccc}
model\hspace{0.5em} & $a$ & $b$ & $\sigma_X$ & $\sigma_Y$\hspace{0.5em} & $\mathbb E[X_t]$ & $\operatorname{Var}(X_t)$ & $\mathcal T_X$\\[2pt]
\hline
multiplicative instability\hspace{0.5em} & $0.65$ & $1.0$ & $0.68$ & $0.55$\hspace{0.5em} & $-0.0004$ & $0.472$ & $1.014$\\
additive forcing\hspace{0.5em} & $1.10$ & $0.516$ & $0.68$ & $0.698$\hspace{0.5em} & $0.003$ & $0.476$ & $1.007$\\
\end{tabular}
\end{ruledtabular}
\end{table}

Throughout Section~\ref{sec:nonlinear}, the model parameters are selected through numerical optimization to reduce discrepancies between selected empirical observable statistics of the two models.
The reported empirical model statistics are computed by averaging the statistics of twenty-four independent long trajectories for each model. The values above show that these parameter choices in this test yield closely matched values of the selected statistics for the two models.

\begin{figure}[tbp!]
\centering
\includegraphics[width=1.0\textwidth]{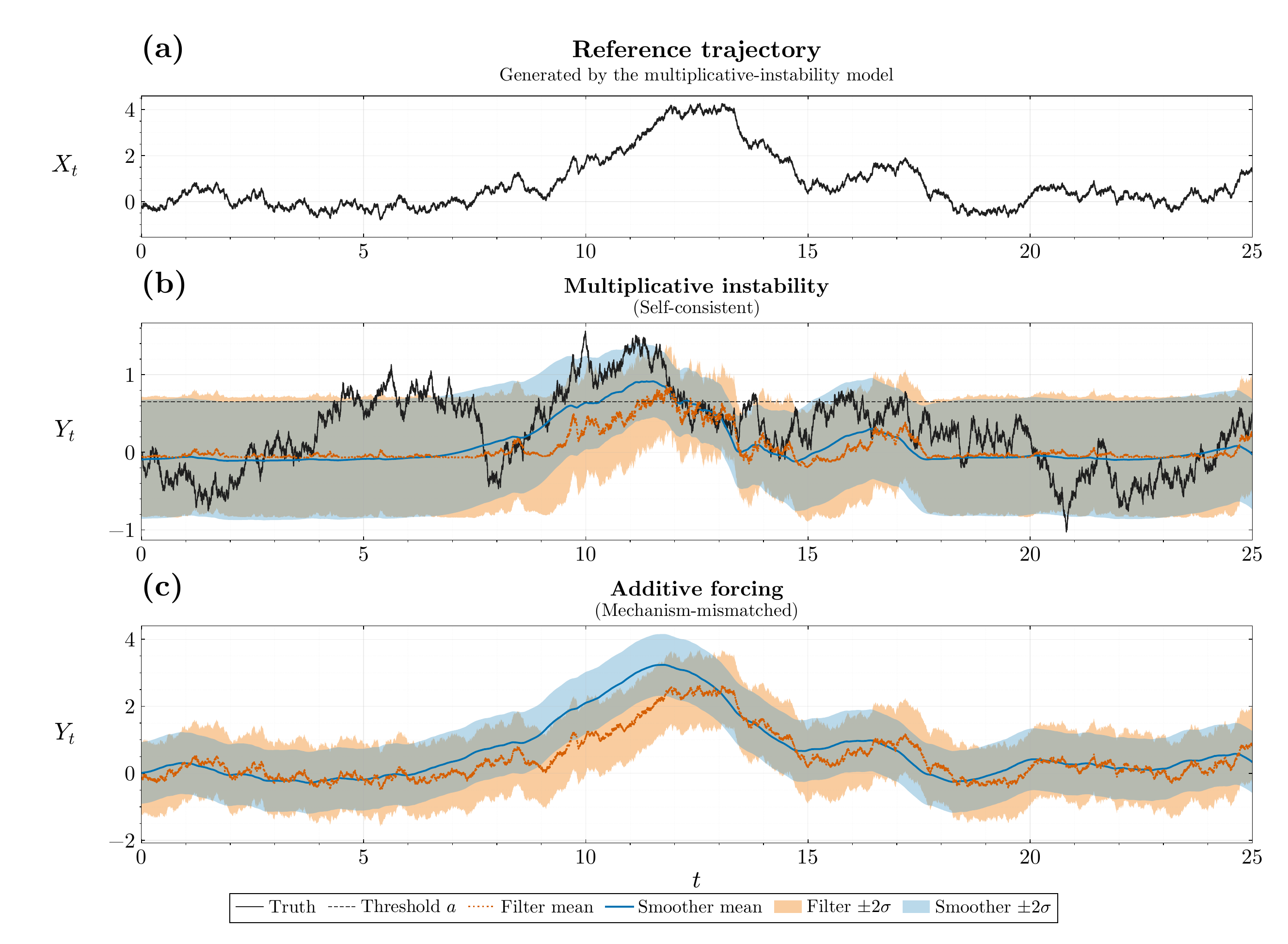}
\caption{Reference trajectory generated by the multiplicative-instability model and hidden-state estimates under both models. Panel (a) shows the \(25\)-time-unit reference trajectory \(X\). Panels (b) and (c) show estimates of \(Y_t\) under the multiplicative-instability model (self-consistent) and the additive-forcing model (mechanism-mismatched), respectively, using the same observed samples. Orange dotted  and blue solid curves show the filtering and 25-time-unit fixed-lag smoothing means, with shaded bands indicating \(\pm2\) posterior standard deviations. In panel (b), the solid black curve shows the true hidden trajectory and the dashed black line marks the instability threshold \(a\).}
\label{fig:sec4_extreme_mult_state_diagnostics}
\end{figure}

\begin{figure}[tbp!]
\centering
\includegraphics[width=1.0\textwidth]{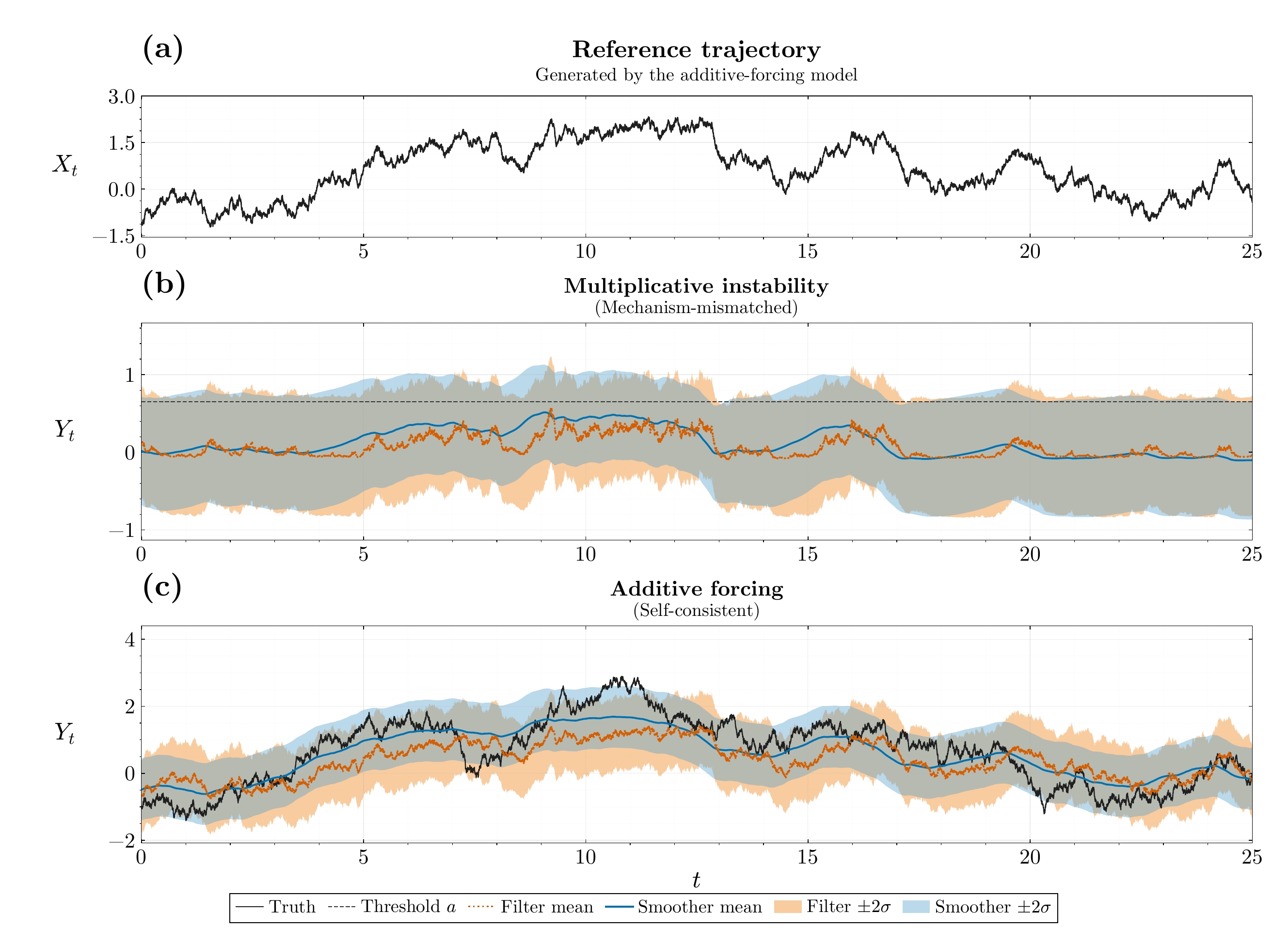}
\caption{Similar to Figure~\ref{fig:sec4_extreme_mult_state_diagnostics} but for a reference trajectory generated by the additive-forcing model. The self-consistent estimates and true hidden trajectory appear in panel (c).}
\label{fig:sec4_extreme_add_state_diagnostics}
\end{figure}

Both models have the CGNS structure, so we compute their filtering and smoothing means and variances using the exact equations in Appendix~\ref{app:cgns_derivations} and substitute them into the Gaussian KL formula \eqref{eq:sec3_gaussian_kl} to obtain the causal fingerprints. Figures~\ref{fig:sec4_extreme_mult_state_diagnostics} and~\ref{fig:sec4_extreme_add_state_diagnostics} show the \(25\)-time-unit reference trajectories generated by the multiplicative-instability and additive-forcing models, respectively, together with the hidden-state filtering and smoothing estimates obtained under both models.

\begin{figure}[tbp!]
\centering
\includegraphics[width=1.0\textwidth]{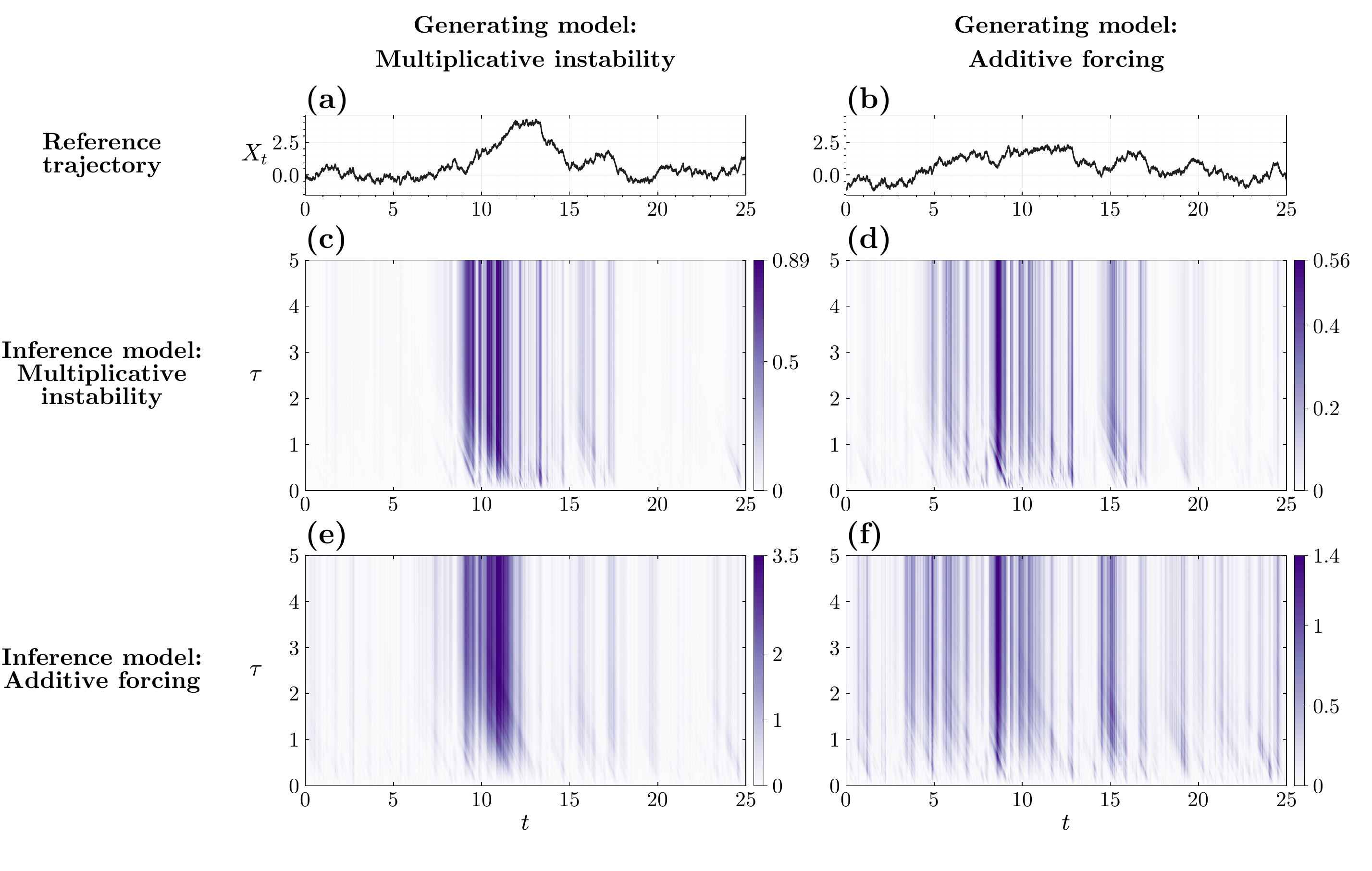}
\caption{Causal fingerprints \(\mathcal I(t,\tau)\) for the extreme-event test. Column headings identify generating models (multiplicative instability, left; additive forcing, right); labels on the left identify inference models for the two heat-map rows in the same order. Panels (a) and (b) show the reference trajectories from Figures~\ref{fig:sec4_extreme_mult_state_diagnostics} and~\ref{fig:sec4_extreme_add_state_diagnostics}, respectively; both fingerprint heat maps in each column are computed from the reference trajectory shown above. Panels (c) and (f) are self-consistent cases; panels (d) and (e) are mechanism-mismatched cases.}
\label{fig:sec4_extreme_fingerprint_diagnostics}
\end{figure}

\paragraph{Results.}
We first consider the reference trajectory generated by the multiplicative-insta\-bility model, with fingerprints under both inference models shown in panels (c) and (e) of Figure~\ref{fig:sec4_extreme_fingerprint_diagnostics}. In the self-consistent case, the multiplicative-instability model's fingerprint captures the onset and early development of the burst through prominent responses over \(t\!\approx\!9\)\,--\,\(t\!\approx\!11\). During this interval, the true hidden state \(Y_t\) rises above the instability threshold \(a\), driving growth toward the large-amplitude burst around \(t\approx12.5\) (panels (a) and (b) of Figure~\ref{fig:sec4_extreme_mult_state_diagnostics}). In the mechanism-mismatched case, the additive-forcing model's prominent fingerprint responses appear later, around \(t\!\approx\!10\)\,--\,\(t\!\approx\!12\), closer to the observed peak. This timing shift is consistent with the distinct mechanisms through which the models explain the burst: in the multiplicative-instability model, an earlier sustained instability allows exponential amplification to build up the burst over time; in the additive-forcing model, strong hidden forcing directly drives local changes in \(X_t\), with fingerprint responses more closely aligned with the burst's growth. In both cases, the prominent fingerprint responses do not persist through the decay of the burst around \(t\!\approx13\!\)\,--\,\(t\!\approx\!15\). As the burst develops, the observations of \(X_t\) already tightly constrain the filtering estimate during the decay. The smoother therefore adds little information, and the fingerprint remains weak.

Next, we consider the reference trajectory generated by the additive-forcing model, with fingerprints under both inference models, shown in panels (d) and (f) of Figure~\ref{fig:sec4_extreme_fingerprint_diagnostics}. Both models produce prominent fingerprint responses at similar state times, aligned with local changes in the observed trajectory, such as the turning points of \(X_t\) near \(t\approx8.5\) and \(t\approx15\). The additive-forcing model explains these local changes through direct hidden forcing, whereas the multiplicative-instability model attributes them to modulation of the effective damping. Under the multiplicative-instability model, both the filtering and smoothing mean estimates remain below the instability threshold throughout the time window (panel (b) of Figure~\ref{fig:sec4_extreme_add_state_diagnostics}), consistent with an explanation of these local changes through damping modulation without entering the unstable regime.

\subsection{Oscillatory variability: nonlinear feedback versus state-dependent stochastic forcing}
\label{sec:enso}

Oscillatory variability can be strongly asymmetric between its positive and negative phases. Similar observable asymmetry can arise from different latent dynamical mechanisms, including nonlinear deterministic feedback and stochastic forcing whose amplitude depends on the system state. These mechanisms may reproduce comparable low-order statistics while assigning different roles to hidden variables during individual events. They therefore provide a useful test of whether causal fingerprints can separate competing explanations of oscillatory behavior.

As a concrete case study, we consider the El Ni\~no-Southern Oscillation (ENSO) \cite{philander1983nino, clarke2008introduction}, a coupled ocean-atmosphere mode of variability in the tropical Pacific with characteristic periods of roughly \(3\)--\(7\) years~\cite{vialard2025recharge}. Its warm phase, El Ni\~no, and cold phase, La Ni\~na, are associated with opposite sea-surface-temperature (SST) anomalies and changes in upper-ocean heat content. The slow recharge and discharge of equatorial heat content helps organize transitions between these phases~\cite{jin1997recharge1}. ENSO also exhibits pronounced warm-cold asymmetry, including positive SST skewness in the eastern equatorial Pacific, while stochastic atmospheric forcing can modify the development and transition of individual events~\cite{vialard2025recharge}.

The {recharge oscillator} (RO) provides a low-dimensional framework for representing these coupled SST and upper-ocean dynamics~\cite{jin2020simpleenso,vialard2025recharge}. Extensions of the RO model can reproduce ENSO asymmetry through different mechanisms, including nonlinear deterministic feedbacks and state-dependent stochastic forcing. We use these two alternatives to examine whether causal fingerprints reveal how competing mechanisms assign hidden influence during the same type of oscillatory event.

Simple RO models typically use an SST anomaly \(T_t\) and an upper-ocean heat-content or thermocline-depth anomaly \(h_t\). The precise definitions of these variables vary across RO formulations~\cite{vialard2025recharge}. Here, \(T_t\) is treated as the observed variable, while \(h_t\) is a hidden recharge variable expressed as a thermocline-depth anomaly. The variables \(T_t\) and \(h_t\) are measured in kelvin and meters, respectively, and time is measured in years.

The first model represents {nonlinear deterministic feedback}, combining linear coupling between \(T_t\) and \(h_t\) with quadratic and cubic SST-anomaly feedback terms:
\begin{equation}
\begin{aligned}
dT_t &=
\left(C_T+RT_t+F_1h_t+b_2T_t^2+b_3T_t^3\right)dt+\sigma_TdW_t^T,\\
dh_t &=
\left(-\epsilon h_t-F_2T_t\right)dt+\sigma_h dW_t^h .
\end{aligned}
\label{eq:sec4_nonlinear_enso}
\end{equation}
The four linear terms form the oscillatory skeleton of the RO model. With \(R<0\), \(RT_t\) represents net linear damping of SST anomalies, weakening both warm and cold anomalies. The \(F_1h_t\) term represents delayed oceanic feedback: positive \(h_t\) promotes SST warming, while negative \(h_t\) favors cooling \cite{jinan1999feedbacks}. The \(-\epsilon h_t\) term relaxes the recharge anomaly toward zero, while \(-F_2T_t\) represents discharge during El Ni\~no and recharge during La Ni\~na, contributing to subsequent phase transitions \cite{jin1997recharge2,vialard2025recharge}.

The quadratic term \(b_2T_t^2\) accounts for nonlinear feedbacks that generate warm-cold asymmetry of SST anomalies, including stronger deep atmospheric convection during El Ni\~no, nonlinear dynamical heating~\cite{jin2003nonlinearheating}, and thermocline-feedback nonlinearity \cite{vialard2025recharge}. With \(b_2>0\), the quadratic term breaks warm-cold symmetry by adding a warming tendency for both positive and negative SST anomalies. The cubic term with \(b_3<0\) provides nonlinear saturation at large amplitudes. The constant drift \(C_T\) is chosen to offset the nonzero mean induced by the nonlinear feedbacks and maintain the zero-mean SST-anomaly convention. The independent additive noise terms \(\sigma_T\,dW_t^T\) and \(\sigma_h\,dW_t^h\) represent unresolved fast forcing, including short-term wind-stress and heat-flux perturbations \cite{vialard2025recharge}.

The second model represents {linear drift with multiplicative noise}. It retains the linear recharge-discharge skeleton and introduces state-dependent noise in the SST equation:
\begin{equation}
\begin{aligned}
dT_t &= (RT_t+F_1h_t)dt+\sigma_T(T_t)dW_t^T,\\
dh_t &= (-\epsilon h_t-F_2T_t)dt+\sigma_h dW_t^h,
\end{aligned}
\label{eq:sec4_linear_mult_noise}
\end{equation}
where
\[
\sigma_T(T)
=\sigma_0\left(1+B\max(T,0)\right).
\]
The state-dependent multiplicative noise \(\sigma_T(T_t)dW_t^T\) accounts for fast atmospheric forcing associated with {westerly wind bursts}, where brief episodes of anomalous westerly winds can last days up to weeks~\cite{levinejin2017noise}. These atmospheric-wind feedbacks tend to be more frequent and stronger during El Ni\~no, strengthening stochastic forcing and contributing to warm-cold asymmetry~\cite{vialard2025recharge}.
Here, \(\sigma_0\) sets the baseline SST-noise amplitude for \(T_t\le0\), while \(B>0\) controls the warm-state enhancement of the noise amplitude.

\paragraph{Parameter choices and experimental setup.}
Parameters of the two models are chosen to reproduce selected observational statistics estimated from the observational SST data, namely the mean, variance, and decorrelation time. These statistics are
\[
\mathbb E[T_t]=0\,\mathrm{K},\qquad
V_T=0.8983\,\mathrm{K}^2,\qquad
\mathcal T_T=0.4682\,\mathrm{yr},
\]
obtained from monthly SST anomalies based on the NCEP Global Ocean Data Assimilation System (GODAS) reanalysis dataset~\cite{behringer2004evaluation}, which are spatially-averaged over the Ni\~no-3 region (eastern equatorial Pacific: \(150^\circ\mathrm{W}\)--\(90^\circ\mathrm{W}\), \(5^\circ\mathrm{S}\)--\(5^\circ\mathrm{N}\)) from January 1980 through March 2026. The corresponding parameters of the two models are listed in Table~\ref{tab:enso_parameters}.
\begin{table}[H]
\centering
\caption{Parameters of the two ENSO models. A dash indicates a parameter not used in the corresponding model.}
\label{tab:enso_parameters}
\renewcommand{\arraystretch}{1.15}
\begin{ruledtabular}
\begin{tabular}{cccc}
parameter & unit\hspace{0.5em} & nonlinear deterministic feedback & linear drift with multiplicative noise\\[2pt]
\hline
$R$ & $\mathrm{yr}^{-1}$\hspace{0.5em} & $-0.839$ & $-0.730$\\
$\epsilon$ & $\mathrm{yr}^{-1}$\hspace{0.5em} & $0.38$ & $0.38$\\
$F_1$ & $\mathrm{K}\,\mathrm{m}^{-1}\mathrm{yr}^{-1}$\hspace{0.5em} & $0.19$ & $0.19$\\
$F_2$ & $\mathrm{m}\,\mathrm{K}^{-1}\mathrm{yr}^{-1}$\hspace{0.5em} & $17.58$ & $13.51$\\
$\sigma_h$ & $\mathrm{m}\,\mathrm{yr}^{-1/2}$\hspace{0.5em} & $6.614$ & $5.942$\\
$C_T$ & $\mathrm{K}\,\mathrm{yr}^{-1}$\hspace{0.5em} & $-0.895$ & ---\\
$b_2$ & $\mathrm{K}^{-1}\mathrm{yr}^{-1}$\hspace{0.5em} & $1.210$ & ---\\
$b_3$ & $\mathrm{K}^{-2}\mathrm{yr}^{-1}$\hspace{0.5em} & $-0.176$ & ---\\
$\sigma_T$ & $\mathrm{K}\,\mathrm{yr}^{-1/2}$\hspace{0.5em} & $1.022$ & ---\\
$\sigma_0$ & $\mathrm{K}\,\mathrm{yr}^{-1/2}$\hspace{0.5em} & --- & $1.028$\\
$B$ & $\mathrm{K}^{-1}$\hspace{0.5em} & --- & $0.392$\\
\end{tabular}
\end{ruledtabular}
\end{table}
With these parameter choices, the two models yield closely matched values of these empirical SST-anomaly statistics, as listed in Table~\ref{tab:enso_statistics}.
\begin{table}[H]
\centering
\caption{Empirical SST-anomaly statistics of the two ENSO models.}
\label{tab:enso_statistics}
\renewcommand{\arraystretch}{1.15}
\begin{ruledtabular}
\begin{tabular}{cccc}
model\hspace{0.5em} & $\mathbb E[T_t]\;(\mathrm{K})$ & $V_T\;(\mathrm{K}^2)$ & $\mathcal T_T\;(\mathrm{yr})$\\[2pt]
\hline
nonlinear deterministic feedback\hspace{0.5em} & $0.0003$ & $0.894$ & $0.466$\\
linear drift with multiplicative noise\hspace{0.5em} & $-0.001$ & $0.891$ & $0.471$\\
\end{tabular}
\end{ruledtabular}
\end{table}

\begin{figure}[tbp!]
\centering
\includegraphics[width=1.0\textwidth]{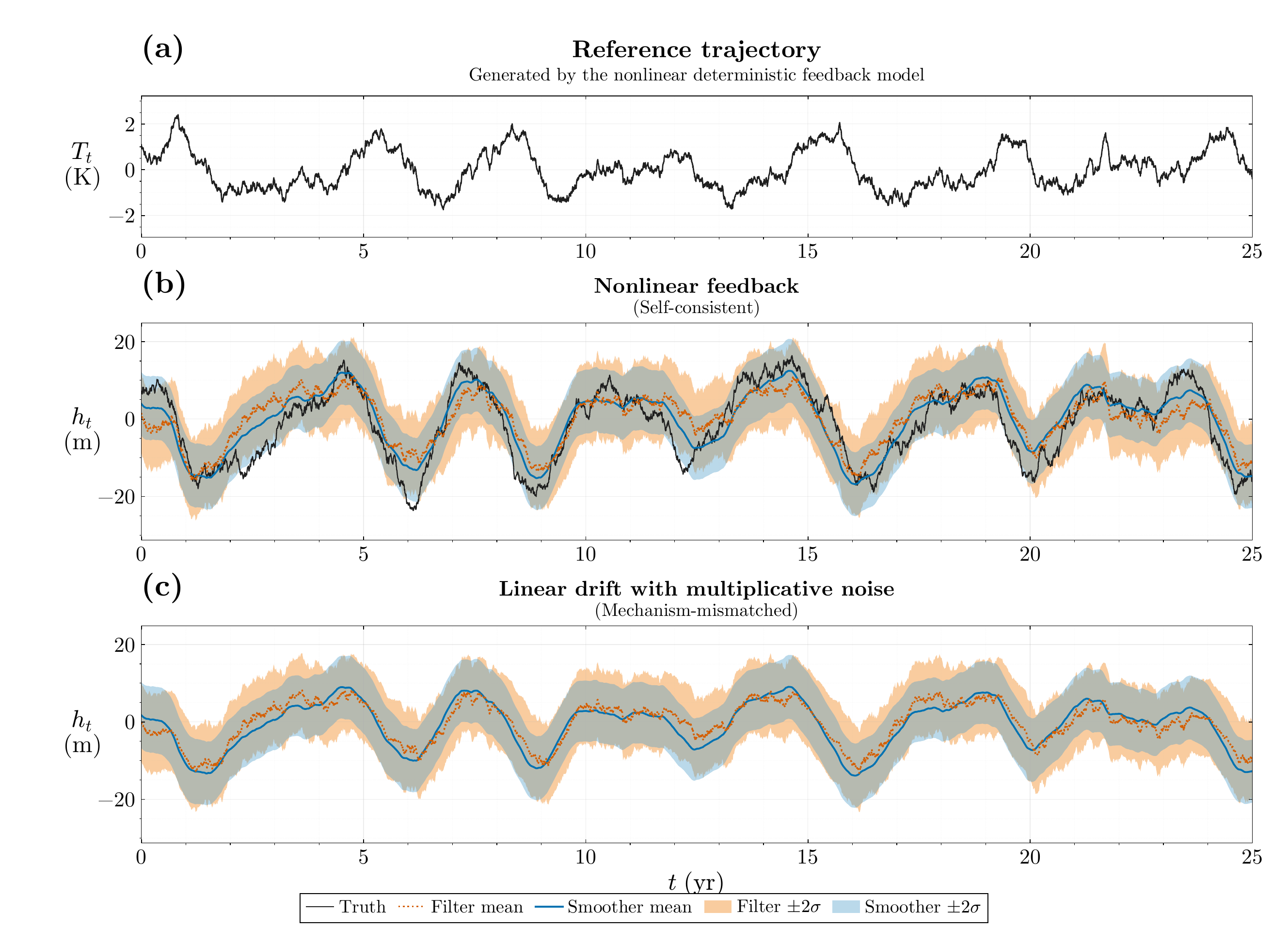}
\caption{Reference trajectory generated by the nonlinear deterministic feedback model and hidden-state estimates under both models. Panel (a) shows the \(25\)-year SST-anomaly reference trajectory \(T\). Panels (b) and (c) show estimates of the recharge anomaly \(h_t\) under the nonlinear deterministic feedback model (self-consistent) and the linear drift with multiplicative noise model (mechanism-mismatched), respectively, using the same SST observations. Orange dotted and blue solid curves show the filtering and 25-year fixed-lag smoothing means, with shaded bands indicating \(\pm2\) posterior standard deviations. The solid black curve in panel (b) shows the true recharge anomaly. \(T\) and \(h\) are measured in kelvin and meters, respectively.}
\label{fig:sec4_enso_nonlinear_state_diagnostics}
\end{figure}

\begin{figure}[tbp!]
\centering
\includegraphics[width=1.0\textwidth]{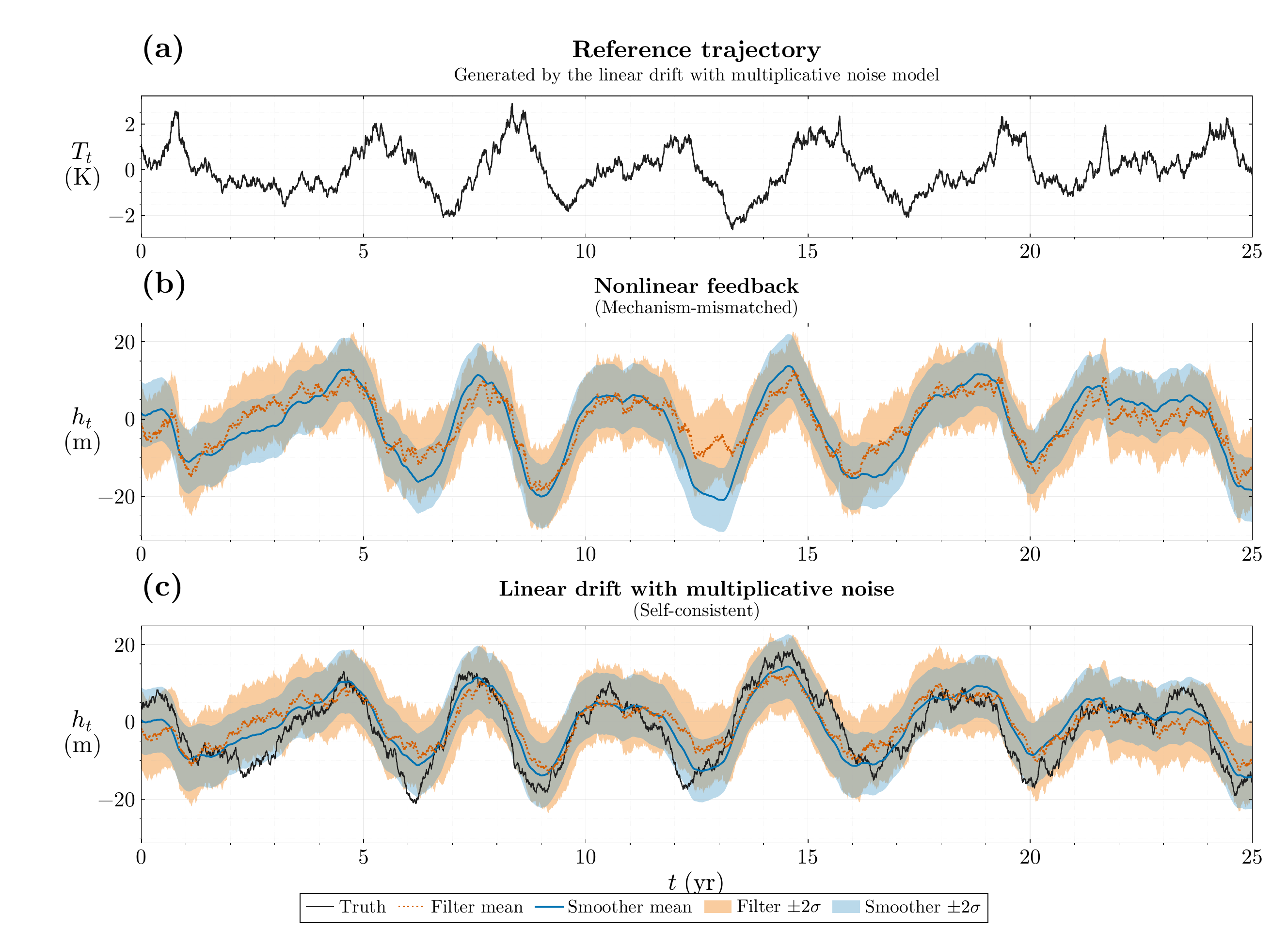}
\caption{Similar to Figure~\ref{fig:sec4_enso_nonlinear_state_diagnostics} but for a reference SST-anomaly trajectory generated by the linear drift with multiplicative noise model. The self-consistent estimates and true recharge trajectory appear in panel (c).}
\label{fig:sec4_enso_mult_noise_state_diagnostics}
\end{figure}

The numerical integration time step is \(\Delta t=5\times10^{-3}\,\mathrm{yr}\). As in Section~\ref{sec:extreme-events}, each dynamical model, \eqref{eq:sec4_nonlinear_enso} and \eqref{eq:sec4_linear_mult_noise}, generates a \(25\)-year reference SST-anomaly trajectory with the same random seed. Then, for each trajectory, fingerprints are computed under both models, yielding self-consistent and mechanism-mismatched results for each one, using the exact CGNS filtering and smoothing equations in Appendix~\ref{app:cgns_derivations} and the Gaussian KL formula \eqref{eq:sec3_gaussian_kl}. Figures~\ref{fig:sec4_enso_nonlinear_state_diagnostics} and~\ref{fig:sec4_enso_mult_noise_state_diagnostics} show these reference trajectories and the corresponding estimates of the recharge anomaly \(h_t\) under both models.

\begin{figure}[tbp!]
\centering
\includegraphics[width=1.0\textwidth]{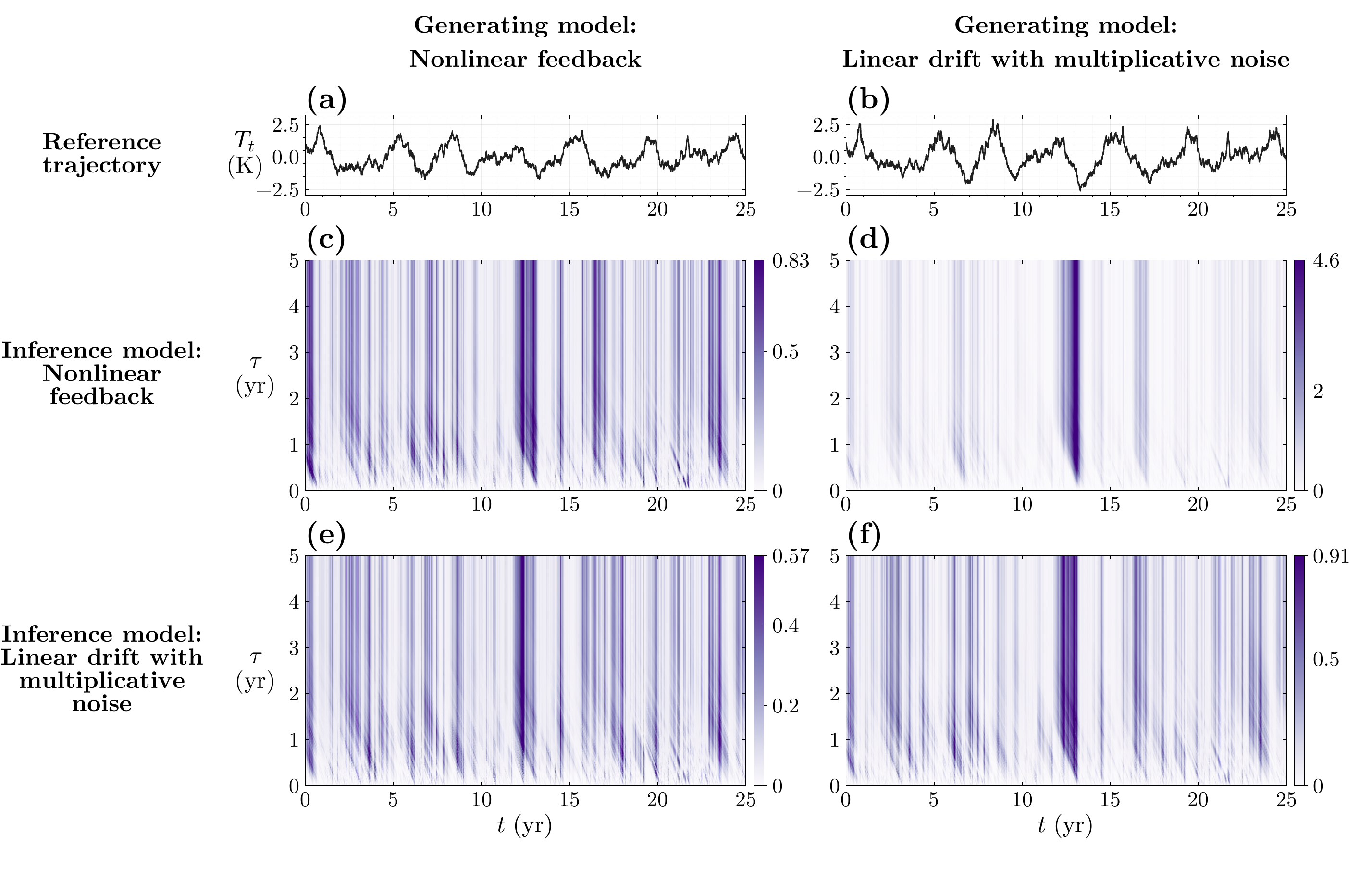}
\caption{Causal fingerprints \(\mathcal I(t,\tau)\) for the ENSO test. Column headings identify generating models (nonlinear deterministic feedback, left; linear drift with multiplicative noise, right); labels on the left identify inference models for the two heat-map rows in the same order. Panels (a) and (b) show the reference SST-anomaly trajectories from Figures~\ref{fig:sec4_enso_nonlinear_state_diagnostics} and~\ref{fig:sec4_enso_mult_noise_state_diagnostics}, respectively; both fingerprint heat maps in each column are computed from the reference trajectory shown above. Panels (c) and (f) are self-consistent cases; panels (d) and (e) are mechanism-mismatched cases.}
\label{fig:sec4_enso_heatmaps}
\end{figure}

\paragraph{Results.}
In the self-consistent cases, the fingerprint responses occur at similar state times \(t\) (panels (c) and (f) of Figure~\ref{fig:sec4_enso_heatmaps}). In both self-consistent cases, prominent fingerprint responses occur near local minima and maxima of the inferred recharge anomaly \(h_t\), where strong oceanic cooling and warming support La Ni\~na and El Ni\~no development, respectively. Examples are near \(t\!\approx\!12\,\mathrm{yr}\)\,--\,\(t\!\approx\!13\,\mathrm{yr}\) and \(t\!\approx\!14\,\mathrm{yr}\)\,--\,\(t\!\approx\!15\,\mathrm{yr}\) in panel (b) of Figure~\ref{fig:sec4_enso_nonlinear_state_diagnostics}. For the reference trajectory generated by the nonlinear deterministic feedback model, the linear drift with multiplicative noise model (panel (e) of Figure~\ref{fig:sec4_enso_heatmaps}; mechanism-mismatched) also produces prominent fingerprint responses at state times similar to those in the self-consistent case (panel (c) of Figure~\ref{fig:sec4_enso_heatmaps}).

A markedly different pattern emerges when the generating model is the one with linear drift with multiplicative noise and the inference model is the one with nonlinear deterministic feedback (panel (d) of Figure~\ref{fig:sec4_enso_heatmaps}; mechanism-mismatched). The fingerprint is dominated by a strong band near \(t\!\approx\!12\,\mathrm{yr}\)\,--\,\(t\!\approx\!13\,\mathrm{yr}\), with weaker responses near \(t\!\approx\!6\,\mathrm{yr}\)\,--\,\(t\!\approx\!7\,\mathrm{yr}\) and \(t\!\approx\!16\,\mathrm{yr}\)\,--\,\(t\!\approx\!17\,\mathrm{yr}\). This band aligns with the development of a strong La Ni\~na event whose SST anomaly reaches its minimum near \(t\approx13\,\mathrm{yr}\) in panel (a) of Figure~\ref{fig:sec4_enso_mult_noise_state_diagnostics}.

Both deterministic nonlinear feedback and state-dependent stochastic forcing have been used in practice to account for positive SST skewness, an essential feature of the equatorial eastern Pacific region~\cite{vialard2025recharge}. The two mechanisms, however, act differently during La Ni\~na development. The quadratic term \(b_2T^2_t\) adds a positive tendency for either sign of \(T_t\), amplifying warm excursions while opposing strong cold excursions and thereby contributing to positive SST skewness. In contrast, multiplicative noise enhances stochastic forcing during El Ni\~no and remains at its baseline amplitude during La Ni\~na. Under the nonlinear-feedback inference model, the enhanced fingerprint during La Ni\~na development is associated with a stronger inferred cooling contribution through \(F_1h_t\), which offsets the quadratic warming tendency \(b_2T_t^2\). Previous RO studies have found a larger contribution to SST skewness from deterministic nonlinear feedbacks than from state-dependent stochastic forcing \cite{vialard2025recharge}.  The fingerprint therefore exposes how the two models assign different hidden-state contributions during La Ni\~na development.

\subsection{Regime switching: deterministic bistability versus hidden double-well forcing}\label{sec:regime-switching}

We finally consider regime-switching systems, in which the observed variable remains near one of two persistent states for extended periods and occasionally transitions between them \cite{majda2012lessons}. Such switching can arise through two distinct mechanisms: {intrinsic bistability}, in which the nonlinear dynamics of the observed variable \(X_t\) itself form a double well and the hidden stochastic forcing \(Y_t\) can trigger transitions of \(X_t\) between the two wells~\cite{hanggi1990reaction}, or {hidden double-well forcing}, in which the dynamics of \(X_t\) are linear and the double-well structure resides in the hidden forcing instead.

The deterministic-bistability model places the double-well structure in the deterministic dynamics of the observed variable:
\begin{equation}
\begin{aligned}
dX_t &= (aX_t-\rho X_t^3+Y_t)dt + \sigma_XdW_t^X, \\
dY_t &= -bY_t\,dt + \sigma_YdW_t^Y .
\end{aligned}
\label{eq:sec4_bistable}
\end{equation}
When \(Y_t=0\), the drift \(aX_t-\rho X_t^3\) has an unstable equilibrium at \(X=0\) and two stable equilibria at \(X=\pm\sqrt{a/\rho}\). The hidden OU process \(Y\) acts as a stochastic forcing: a sufficiently strong fluctuation in \(Y_t\) can push \(X_t\) across the unstable equilibrium into the opposite basin, where the bistable drift stabilizes \(X_t\) in the new regime.

In contrast, the hidden-double-well-forcing model treats \(Z_t\) as a hidden stochastic forcing with its own double-well dynamics:
\begin{equation}
\begin{aligned}
dX_t &= (-\alpha X_t+Z_t)dt + \sigma_XdW_t^X, \\
dZ_t &= (a_Z Z_t-\rho_Z Z_t^3)dt + \sigma_ZdW_t^Z .
\end{aligned}
\label{eq:sec4_bistable_forcing}
\end{equation}
When \(a_Z,\rho_Z>0\), the drift \(a_Z Z_t-\rho_Z Z_t^3\) has two stable equilibria at \(Z=\pm\sqrt{a_Z/\rho_Z}\). For a fixed \(Z_t\), the dynamics of \(X_t\) are linear and stable, relaxing toward \(Z_t/\alpha\) with no double-well structure of their own. When \(Z_t\) switches to the opposite well, its forcing drives \(X_t\) toward the corresponding observed regime. As long as \(Z_t\) remains in that well, its forcing sustains the new observed regime.

\paragraph{Parameter choices and experimental setup.}
The parameters for the two models are
\[
\begin{aligned}
\text{deterministic bistability:}\quad
&a=0.85,\quad \rho=0.3778,\quad b=0.40,\quad \sigma_X=0.25,\quad \sigma_Y=0.35,\\
\text{hidden double-well forcing:}\quad
&\alpha=0.7850,\quad a_Z=1.140,\quad \rho_Z=0.7401,\quad \sigma_X=0.25,\quad \sigma_Z=0.5304.
\end{aligned}
\]
With these parameters, the deterministic-bistability model has equilibria at \(X=\pm1.5\). The hidden-double-well-forcing model has equilibria at \(Z\approx\pm1.241\), which induce observed equilibria at \(X\approx\pm1.581\).

These parameter choices yield closely-matched empirical stationary means, variances and decorrelation times for the observable \(X\). The corresponding empirical equilibrium statistics of the mean, variance, and decorrelation time are listed in Table~\ref{tab:bistability_statistics}.
\begin{table}[H]
\centering
\caption{Empirical observable statistics of the regime-switching models.}
\label{tab:bistability_statistics}
\renewcommand{\arraystretch}{1.15}
\begin{ruledtabular}
\begin{tabular}{cccc}
model\hspace{0.5em} & $\mathbb E[X_t]$ & $\operatorname{Var}(X_t)$ & $\mathcal T_X$\\[2pt]
\hline
deterministic bistability\hspace{0.5em} & $0.07768$ & $2.1072$ & $9.3554$\\
hidden double-well forcing\hspace{0.5em} & $0.05445$ & $2.1066$ & $9.3624$\\
\end{tabular}
\end{ruledtabular}
\end{table}

\begin{figure}[tbp!]
\centering
\includegraphics[width=1.0\textwidth]{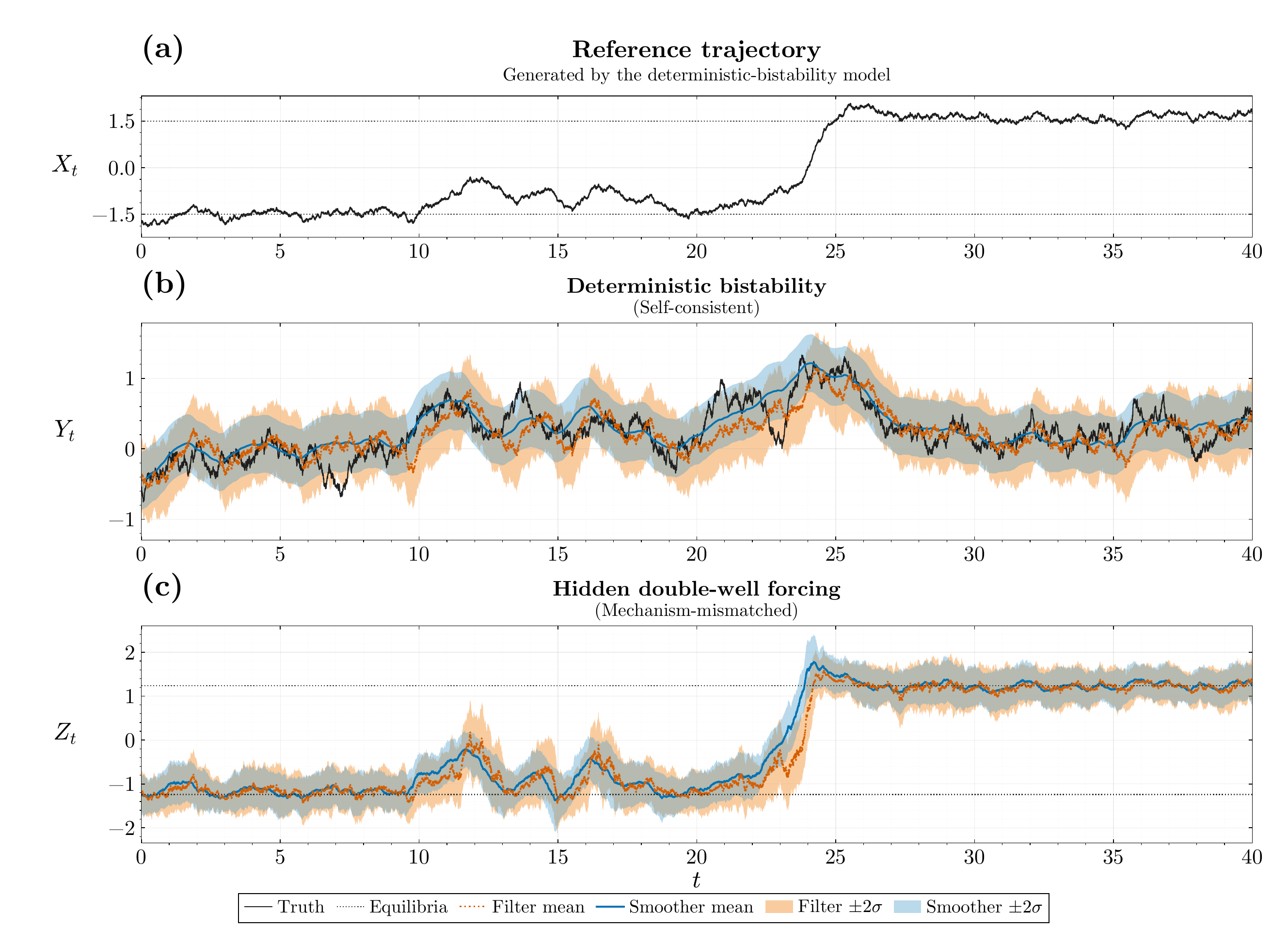}
\caption{Reference trajectory generated by the deterministic-bistability model and hidden-state estimates under both models. Panel (a) shows the \(40\)-time-unit reference trajectory \(X\), with dotted lines marking its stable equilibria under the generating model. Panel (b) shows estimates of \(Y_t\) under the deterministic-bistability model (self-consistent); panel (c) shows estimates of \(Z_t\) under the hidden-double-well-forcing model (mechanism-mismatched), using the same observed samples. Orange dotted and blue solid curves show the filtering and 40-time-unit fixed-lag smoothing means, with shaded bands indicating \(\pm2\) posterior standard deviations. The solid black curve in panel (b) shows the true hidden trajectory. Dotted lines in panel (c) mark the hidden double-well equilibria.}
\label{fig:sec4_bistability_det_state_diagnostics}
\end{figure}

\begin{figure}[tbp!]
\centering
\includegraphics[width=1.0\textwidth]{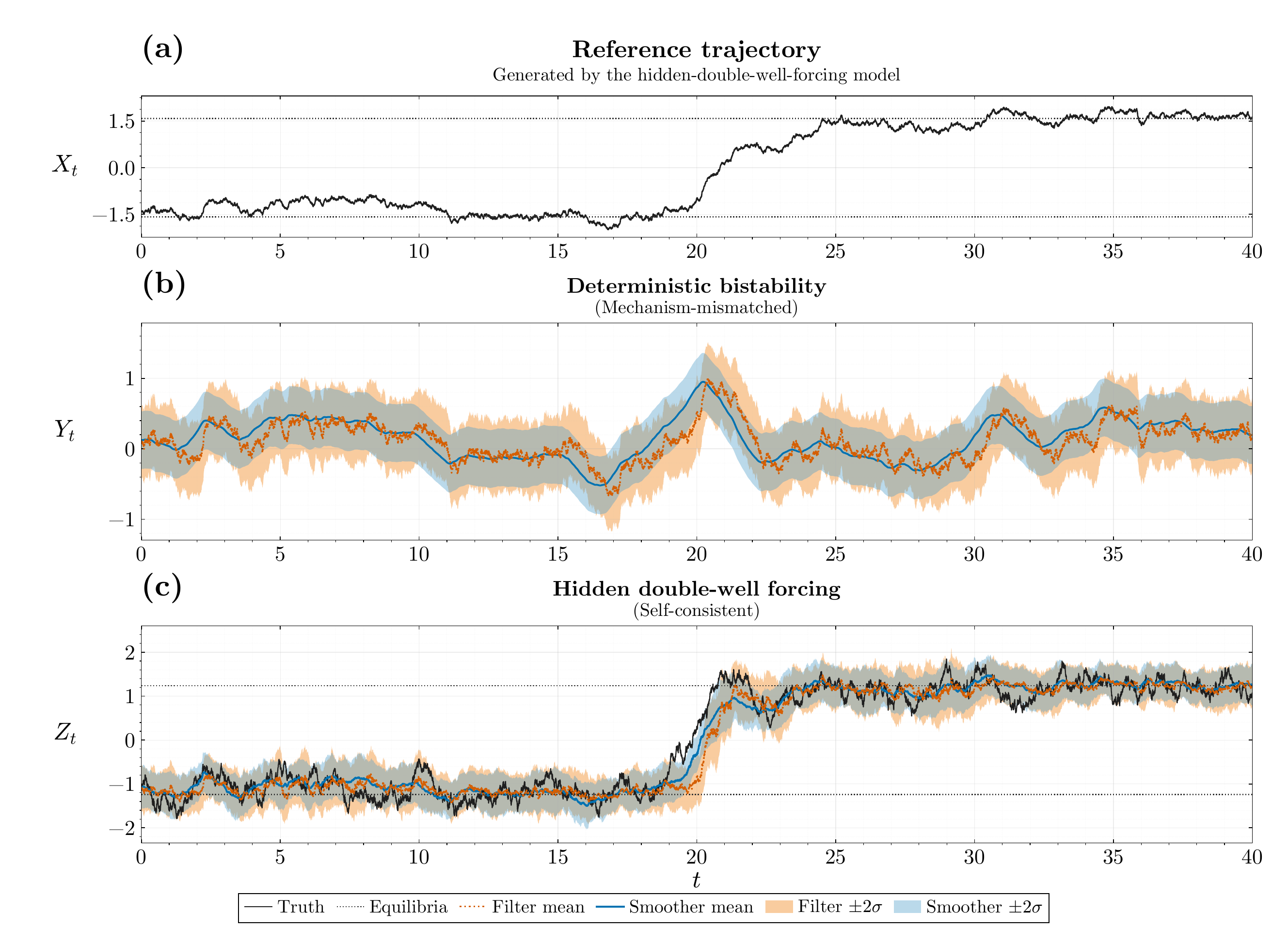}
\caption{Similar to Figure~\ref{fig:sec4_bistability_det_state_diagnostics} but for a reference trajectory generated by the hidden-double-well-forcing model. The self-consistent estimates and true hidden trajectory appear in panel (c); dotted lines mark the corresponding observable and hidden equilibria.}
\label{fig:sec4_bistability_switch_state_diagnostics}
\end{figure}

Each model generates a \(40\)-time-unit reference trajectory \(X\). For each corresponding trajectory, fingerprints are computed under both models. Figures~\ref{fig:sec4_bistability_det_state_diagnostics} and~\ref{fig:sec4_bistability_switch_state_diagnostics} show these reference trajectories together with estimates of \(Y_t\) under the deterministic-bistability model and \(Z_t\) under the hidden-double-well-forcing model.

All numerical simulations in this comparison use the Euler-Maruyama scheme with time step \(\Delta t=2\times10^{-3}\). The deterministic-bistability model has the CGNS structure, so its hidden-state filtering and smoothing distributions are computed using the exact equations in Appendix~\ref{app:cgns_derivations}, while the causal fingerprint is computed using the KL formula in \eqref{eq:sec3_gaussian_kl}. The hidden-double-well-forcing model is not conditionally Gaussian because of its cubic hidden-state drift, so its hidden-state inference uses a \(50\)-member {ensemble Kalman filter}~\cite{evensen1994sequential,burgers1998analysis} and {ensemble Rauch-Tung-Striebel smoother}~\cite{raanes2016ensemble}. For this model, the filtering and smoothing distributions are approximated as Gaussians using their respective ensemble means and variances, and the KL fingerprint is similarly evaluated from \eqref{eq:sec3_gaussian_kl}. This ensemble-Gaussian calculation is an approximation to the theoretical fingerprint in \eqref{eq:info_gain_natural}; in particular, the exact conditional-mutual-information identity \eqref{eq:expected_fingerprint_cmi} applies to exact filtering and smoothing posteriors. We therefore interpret the corresponding fingerprints as approximate diagnostics under the candidate model.

\begin{figure}[tbp!]
\centering
\includegraphics[width=1.0\textwidth]{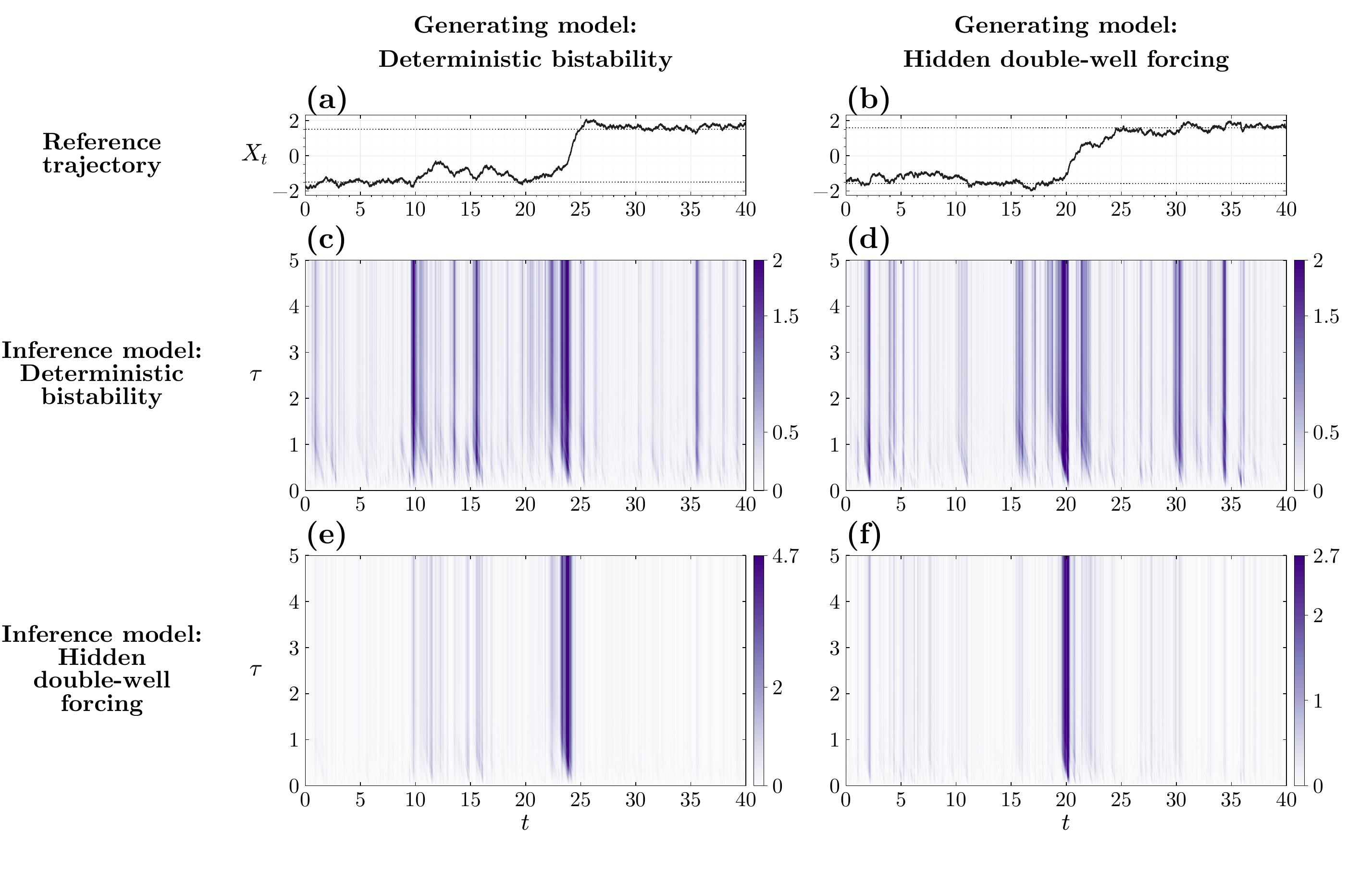}
\caption{Causal fingerprints \(\mathcal I(t,\tau)\) for the regime-switching test. Column headings identify generating models (deterministic bistability, left; hidden double-well forcing, right); labels on the left identify inference models for the two heat-map rows in the same order. Panels (a) and (b) show the reference trajectories from Figures~\ref{fig:sec4_bistability_det_state_diagnostics} and~\ref{fig:sec4_bistability_switch_state_diagnostics}, respectively; both fingerprint heat maps in each column are computed from the reference trajectory shown above. Panels (c) and (f) are self-consistent cases; panels (d) and (e) are mechanism-mismatched cases. Dotted lines mark observable equilibria.}
\label{fig:sec4_bistability_heatmaps}
\end{figure}

\paragraph{Results.}
We first consider the reference trajectory generated by the deterministic-bista\-bility model, with fingerprints under both inference models shown in panels (c) and (e) of Figure~\ref{fig:sec4_bistability_heatmaps}. The reference trajectory in panel (a) of Figure~\ref{fig:sec4_bistability_det_state_diagnostics} undergoes its main transition from the negative to the positive regime, with the tipping point occurring near \(t\approx24\). In the self-consistent case, the deterministic-bistability fingerprint responds near the main transition and also to smaller fluctuations of \(X_t\) within the same regime. In the mechanism-mismatched case, the hidden-double-well-forcing fingerprint is more concentrated near the transition, with relatively much weaker responses to fluctuations within the initial regime. This contrast is consistent with the different hidden dynamics: fluctuations in the OU process \(Y_t\) can trigger regime transitions or intermittent perturbations of \(X_t\) within its current regime, whereas the double-well restoring drift tends to constrain \(Z_t\) near one of its stable equilibria. These relatively weaker within-regime responses indicate that the hidden-double-well-forcing model attributes a smaller causal contribution to hidden forcing during within-regime fluctuations.

The same contrast appears for the reference trajectory generated by the hidden-double-well-forcing model, with the tipping point occurring near \(t\approx20\) (panel (a) of Figure~\ref{fig:sec4_bistability_switch_state_diagnostics}). In panels (d) and (f) of Figure~\ref{fig:sec4_bistability_heatmaps}, the hidden-double-well-forcing fingerprint (self-consistent) is concentrated near the transition, while the deterministic-bistability fingerprint (mechanism-mismatched) also responds to fluctuations within the observed regimes. Thus, causal fingerprints reveal differences between the two regime-switching mechanisms through the contrasting temporal patterns of hidden influence during transitions and fluctuations within a regime.

\section{Discussion and Conclusions}
\label{sec:discussion}

This work introduces an observation-conditioned, model-informed causal fingerprint framework for diagnosing hidden influence in partially-observed dynamical systems. At each state time, the fingerprint measures how future observations revise the inferred hidden state. Evaluated across different future observation windows, it tracks the temporal reach of that influence, while computed along the trajectory, it shows how the pattern changes from one event or event stage to another. Because the observed trajectory is fixed while hidden-state inference depends on the candidate model, fingerprints provide a direct way to compare competing mechanistic explanations of the same observations.

In linear Gaussian systems, the analysis yields an explicit characteristic fingerprint timescale for the temporal reach of hidden causal influence. In Test I, the surrogates reproduce the perfect model's selected observable summary statistics, but removing feedback-induced state memory or suppressing hidden-noise memory changes the strength and temporal structure of their fingerprints. In Test II, the chosen one-component surrogate reproduces features of the fast marginal fingerprint but misses the slow component's contribution over longer future observation windows. These tests demonstrate how fingerprint patterns reveal omitted or misrepresented hidden mechanisms.

The nonlinear tests demonstrate event-level diagnosis of competing hidden mechanisms. For the observed amplitude burst generated by the multiplicative-instability model, fingerprint responses appear earlier under that model, whereas the additive-forcing model's responses are more closely aligned with the burst's growth. In the oscillatory variability case study, distinct responses during La Ni\~na development show how deterministic nonlinear feedback and state-dependent stochastic forcing assign different roles to the hidden recharge state. For regime switching, deterministic-bistability fingerprints respond during transitions and within-regime fluctuations, whereas hidden-double-well-forcing fingerprints are more concentrated near transitions. Across these examples, inference models with similar observable behavior can assign hidden influence to different times and stages of an event.

The fingerprint should be interpreted as a model-relative diagnostic rather than a model-free interventional proof of causality. Its mechanistic interpretation depends on the candidate dynamical structure, the chosen hidden-state representation, and the accuracy of the filtering and smoothing distributions. A fingerprint therefore does not by itself establish that one candidate mechanism is the true physical mechanism; instead, it reveals how competing models assign hidden influence to the same observed evolution. In non-Gaussian settings, approximate posterior methods such as the ensemble calculation in Section~\ref{sec:regime-switching} introduce an additional layer of approximation.

Future work could extend the framework to higher-dimensional nonlinear systems with multiple interacting hidden and observed variables. A key question is how to resolve the joint and individual contributions of hidden processes to different observed variables and trace these influences over time. Another direction is to apply the framework to real-world observations, extending beyond the present ENSO experiment, which uses observational statistics for parameter selection but model-generated trajectories for causal diagnosis. Fingerprints computed along observed SST-anomaly trajectories could help diagnose how alternative mechanisms explain the evolution of individual El Ni\~no and La Ni\~na events. Such analyses could clarify the roles of different physical processes and inform the construction of conceptual and reduced models.

\begin{acknowledgments}
N.C. is supported by the Army Research Office under Award No. W911NF-23-1-0118. Z.J. and M.A. are partially supported as PhD research assistants under this grant. M.A. is also supported in part by NSF Award DMS-2023239.
\end{acknowledgments}

\section*{Author Declarations}
\subsection*{Conflict of Interest}
The authors have no conflicts to disclose.

\subsection*{Author Contributions}
\noindent\textbf{Zhang Jiang}: Formal analysis; Investigation; Methodology; Software; Visualization; Writing -- original draft.

\noindent\textbf{Marios Andreou}: Data Curation; Formal analysis; Methodology; Writing -- review \& editing.

\noindent\textbf{Nan Chen}: Conceptualization; Funding acquisition; Methodology; Project administration; Supervision; Writing -- review \& editing.

\subsection*{AI Statement}
The authors used OpenAI's GPT-5.6 Sol to assist with the language editing of the manuscript. The authors have reviewed the content of this work and take full responsibility for it.

\section*{Data Availability}
The GODAS reanalysis data used in Section~\ref{sec:enso} are publicly available from NOAA~\cite{behringer2004evaluation}. The numerical simulation outputs and the processed SST-anomaly series that support the findings of this study are available from the corresponding author upon reasonable request.

\appendix

\section{Optimal Filtering and Smoothing for Conditional Gaussian Nonlinear Systems}
\label{app:cgns_derivations}

This appendix summarizes the optimal filtering and smoothing equations for conditional Gaussian nonlinear systems (CGNS). These equations are used for the linear Gaussian models in Section~\ref{sec:linear} and the nonlinear models with CGNS structure in Section~\ref{sec:nonlinear}.

\subsection{General conditional Gaussian nonlinear systems}

A {conditional Gaussian nonlinear system} (CGNS) is a nonlinear stochastic model that admits exact Gaussian filtering and smoothing distributions for its hidden variables \cite{liptser2001statistics,liptser2001applications,chenmajda2020smoothing, andreou2024martingale}. Its joint and marginal distributions can nevertheless be non-Gaussian, allowing it to represent features such as intermittency and extreme events \cite{chenmajda2016filtering,chenmajda2018conditional}. We consider general CGNSs of the following form:
\begin{equation}
\begin{aligned}
\mathrm{d}\mathbf{x}(t)
&=
\left(\boldsymbol{\Lambda}_x(t,\mathbf{x})\mathbf{y}(t)+\mathbf{f}_x(t,\mathbf{x})\right)\mathrm{d} t
+\boldsymbol{\Sigma}_x(t,\mathbf{x})\mathrm{d} \mathbf{W}_x(t),\\
\mathrm{d}\mathbf{y}(t)
&=
\left(\boldsymbol{\Lambda}_y(t,\mathbf{x})\mathbf{y}(t)+\mathbf{f}_y(t,\mathbf{x})\right)\mathrm{d} t
+\boldsymbol{\Sigma}_y(t,\mathbf{x})\mathrm{d} \mathbf{W}_y(t),
\end{aligned}
\label{eq:sec3_cgns_general}
\end{equation}
where \(\mathbf{x}(t)\in\mathbb{R}^{d_x}\) is observed, \(\mathbf{y}(t)\in\mathbb{R}^{d_y}\) is hidden, and \(\mathbf{W}_x,\mathbf{W}_y\) are independent Wiener processes. We denote the noise covariances by
\[
\mathbf{Q}_x(t,\mathbf{x}):=\boldsymbol{\Sigma}_x(t,\mathbf{x})\boldsymbol{\Sigma}_x(t,\mathbf{x})^\top,\qquad
\mathbf{Q}_y(t,\mathbf{x}):=\boldsymbol{\Sigma}_y(t,\mathbf{x})\boldsymbol{\Sigma}_y(t,\mathbf{x})^\top,
\]
which are both assumed to be positive definite. We assume that the initial conditional distribution \(p(\mathbf{y}(0)\mid\mathbf{x}(0))\) is Gaussian. Throughout this appendix, every matrix whose inverse appears in the formulas is assumed to be nonsingular. Both equations of \eqref{eq:sec3_cgns_general} are affine in \(\mathbf{y}\), while their coefficients and the noise amplitudes may depend nonlinearly on \(\mathbf{x}\). The linear systems in Section~\ref{sec:linear} are special cases of \eqref{eq:sec3_cgns_general}.

\subsection{Optimal filtering and smoothing equations}

Under suitable regularity conditions, the hidden state has a Gaussian conditional distribution given the observed history \(\mathbf{x}_{[0,t]}:=\{\mathbf{x}(s):0\le s\le t\}\) \cite{liptser2001statistics,liptser2001applications, andreou2024martingale}:
\begin{equation}
\mathbf{y}(t)\mid\mathbf{x}_{[0,t]}
\sim\mathcal N\!\left(\boldsymbol{\mu}_{\mathrm f}(t),\mathbf{R}_{\mathrm f}(t)\right).
\label{eq:app_cgns_filter_distribution}
\end{equation}
This is the filtering distribution, with conditional mean and covariance
\[
\boldsymbol{\mu}_{\mathrm f}(t):=\mathbb{E}[\mathbf{y}(t)\mid \mathbf{x}_{[0,t]}],\qquad
\mathbf{R}_{\mathrm f}(t):=\operatorname{Cov}(\mathbf{y}(t)\mid \mathbf{x}_{[0,t]}).
\]
These moments completely characterize the Gaussian filtering distribution and evolve according to the exact optimal filtering equations:
\begin{equation}
\begin{aligned}
\mathrm{d}\boldsymbol{\mu}_{\mathrm f}(t)
&=
\left(\boldsymbol{\Lambda}_y\boldsymbol{\mu}_{\mathrm f}+\mathbf{f}_y\right)\mathrm{d} t
+\mathbf{K}_{\mathrm f}\left(\mathrm{d}\mathbf{x}-\left(\boldsymbol{\Lambda}_x\boldsymbol{\mu}_{\mathrm f}+\mathbf{f}_x\right)\mathrm{d} t\right),\\
\mathrm{d}\mathbf{R}_{\mathrm f}(t)
&=
\left(\boldsymbol{\Lambda}_y\mathbf{R}_{\mathrm f}+\mathbf{R}_{\mathrm f}\boldsymbol{\Lambda}_y^\top+\mathbf{Q}_y-\mathbf{K}_{\mathrm f}\mathbf{Q}_x\mathbf{K}_{\mathrm f}^\top\right)\mathrm{d} t,
\end{aligned}
\label{eq:sec3_cgns_filter}
\end{equation}
where the model coefficients are evaluated at \((t,\mathbf{x}(t))\), and
\[
\mathbf{K}_{\mathrm f}(t,\mathbf{x}):=\mathbf{R}_{\mathrm f}(t)\boldsymbol{\Lambda}_x(t,\mathbf{x})^\top \mathbf{Q}_x(t,\mathbf{x})^{-1}
\]
is the {Kalman gain}.

Let \(T>0\) denote the fixed terminal observation time, so that observations are available over \([0,T]\). For \(0\le t\le T\), define
\[
\begin{aligned}
\boldsymbol{\mu}_{\mathrm s}(t):=\mathbb{E}[\mathbf{y}(t)\mid \mathbf{x}_{[0,T]}],\qquad
\mathbf{R}_{\mathrm s}(t):=\operatorname{Cov}(\mathbf{y}(t)\mid \mathbf{x}_{[0,T]}).
\end{aligned}
\]
These moments characterize the Gaussian fixed-interval smoothing distribution of \(\mathbf y(t)\) given the observed path \(\mathbf x_{[0,T]}\).
Once the filtering mean and covariance have been computed over \([0,T]\), we initialize the smoother at \(T\) with terminal conditions:
\[
\boldsymbol{\mu}_{\mathrm s}(T)=\boldsymbol{\mu}_{\mathrm f}(T),
\qquad
\mathbf{R}_{\mathrm s}(T)=\mathbf{R}_{\mathrm f}(T).
\]
The optimal smoothing moments are then obtained by integrating the following equations backward from \(T\) to \(0\)~\cite{liptser2001statistics,liptser2001applications, chenmajda2020smoothing, andreou2024martingale}:
\begin{equation}
\begin{aligned}
\overleftarrow{\mathrm{d}}\!\boldsymbol{\mu}_{\mathrm s}(t)
&=
-\left(\boldsymbol{\Lambda}_y\boldsymbol{\mu}_{\mathrm s}+\mathbf{f}_y-\mathbf{Q}_y\mathbf{R}_{\mathrm f}^{-1}(\boldsymbol{\mu}_{\mathrm f}-\boldsymbol{\mu}_{\mathrm s})\right)\mathrm{d} t,\\
\overleftarrow{\mathrm{d}}\!\mathbf{R}_{\mathrm s}(t)
&=
-\left((\boldsymbol{\Lambda}_y+\mathbf{Q}_y\mathbf{R}_{\mathrm f}^{-1})\mathbf{R}_{\mathrm s}
+\mathbf{R}_{\mathrm s}(\boldsymbol{\Lambda}_y+\mathbf{Q}_y\mathbf{R}_{\mathrm f}^{-1})^\top
-\mathbf{Q}_y\right)\mathrm{d} t.
\end{aligned}
\label{eq:sec3_cgns_smoother}
\end{equation}
Here, \(\overleftarrow{\mathrm{d}}\!(\boldsymbol{\cdot})\) is the formal backward-in-time differential, i.e. \(\overleftarrow{\mathrm{d}}\!(\boldsymbol{\cdot}) = -\mathrm{d}(\boldsymbol{\cdot})\).

An equivalent formulation is given by the {online smoother}. While the fixed-interval smoother formulation \eqref{eq:sec3_cgns_smoother} holds the observation endpoint \(T\) fixed and propagates information from future observations backward to earlier hidden-state estimates, the online smoother instead fixes the state time \(t\) and describes how the posterior distribution of \(\mathbf{y}(t)\) changes as additional future observations are incorporated~\cite{andreou2026adaptive}. Define
\[
\boldsymbol{\mu}_{\mathrm s}(t,\tau):=\mathbb{E}[\mathbf{y}(t)\mid \mathbf{x}_{[0,t+\tau]}],\qquad
\mathbf{R}_{\mathrm s}(t,\tau):=\operatorname{Cov}(\mathbf{y}(t)\mid \mathbf{x}_{[0,t+\tau]}).
\]
In the following equations, coefficients written with a single time argument \(s\) are evaluated at \((s,\mathbf x(s))\). Define the {influence kernel} \(\boldsymbol{\Phi}(t,\tau)\) as the fundamental matrix solution of
\[
\frac{\partial}{\partial\tau}\boldsymbol{\Phi}(t,\tau)
=
\left.\left(\boldsymbol{\Lambda}_y-\mathbf{R}_{\mathrm f}\boldsymbol{\Lambda}_x^\top \mathbf{Q}_x^{-1}\boldsymbol{\Lambda}_x\right)\right|_{t+\tau}
\boldsymbol{\Phi}(t,\tau),
\qquad
\boldsymbol{\Phi}(t,0)=\mathbf{I}.
\]
With \(t\) held fixed, we use \(\mathrm d_\tau\) to denote the Itô differential as \(\tau\) increases; in particular, \(\mathrm{d}_\tau\mathbf{x}(t+\tau)\) denotes the observation increment at the advancing endpoint \(t+\tau\). Starting from
\[
\boldsymbol{\mu}_{\mathrm s}(t,0)=\boldsymbol{\mu}_{\mathrm f}(t),
\qquad
\mathbf{R}_{\mathrm s}(t,0)=\mathbf{R}_{\mathrm f}(t),
\]
the online smoother equations are
\begin{equation}
\begin{aligned}
\mathrm{d}_\tau\boldsymbol{\mu}_{\mathrm s}(t,\tau)
&=
\mathbf{R}_{\mathrm f}(t)\boldsymbol{\Phi}(t,\tau)^\top\boldsymbol{\Lambda}_x(t+\tau)^\top \mathbf{Q}_x(t+\tau)^{-1}
\left[
\mathrm{d}_\tau\mathbf{x}(t+\tau)
-\left.\left(\boldsymbol{\Lambda}_x\boldsymbol{\mu}_{\mathrm f}+\mathbf{f}_x\right)\right|_{t+\tau}\mathrm{d}\tau
\right],\\
\mathrm{d}_\tau\mathbf{R}_{\mathrm s}(t,\tau)
&=
-\mathbf{R}_{\mathrm f}(t)\boldsymbol{\Phi}(t,\tau)^\top\boldsymbol{\Lambda}_x(t+\tau)^\top
\mathbf{Q}_x(t+\tau)^{-1}\boldsymbol{\Lambda}_x(t+\tau)\boldsymbol{\Phi}(t,\tau)\mathbf{R}_{\mathrm f}(t)\mathrm{d}\tau.
\end{aligned}
\label{eq:sec3_cgns_online_smoother}
\end{equation}
For \(T=t+\tau\), the fixed-interval and online formulations yield the same smoothing distribution for \(\mathbf{y}(t)\).

\paragraph{Two-dimensional linear Gaussian model.}
For the model \eqref{eq:sec3_perfect_model}, take \(\mathbf{x}=X\), \(\mathbf{y}=Y\), and
\[
\begin{aligned}
\Lambda_x&=1,&\qquad f_x&=-aX,&\qquad Q_x&=\sigma_X^2,\\
\Lambda_y&=-b,&\qquad f_y&=cX,&\qquad Q_y&=\sigma_Y^2.
\end{aligned}
\]
These substitutions give the filtering variance and online smoothing equations in Section~\ref{sec:sec3_closed_form}. The filtering mean and variance equations obtained from \eqref{eq:sec3_cgns_filter} are
\begin{equation}
\begin{aligned}
\mathrm{d}\mu_{\mathrm f}(t)
&=
\left(cX_t-b\mu_{\mathrm f}(t)\right)\mathrm{d} t
+
\frac{R_{\mathrm f}(t)}{\sigma_X^2}
\left[
\mathrm{d} X_t
-
\left(\mu_{\mathrm f}(t)-aX_t\right)\mathrm{d} t
\right],\\
\mathrm{d} R_{\mathrm f}(t)
&=
\left[
-2bR_{\mathrm f}(t)
+\sigma_Y^2
-\frac{R_{\mathrm f}(t)^2}{\sigma_X^2}
\right]\mathrm{d} t.
\end{aligned}
\label{eq:sec3_scalar_filter_mean}
\end{equation}

For the terminal observation time \(T\) defined above, we write
\[
\mu_{\mathrm s}(t)
:=
\mu_{\mathrm s}(t,T-t),
\qquad
R_{\mathrm s}(t)
:=
R_{\mathrm s}(t,T-t).
\]
With terminal conditions
\(\mu_{\mathrm s}(T)=\mu_{\mathrm f}(T)\) and
\(R_{\mathrm s}(T)=R_{\mathrm f}(T)\), the fixed-interval smoother integrates backward in time from \(T\) to \(t\) according to
\begin{equation}
\begin{aligned}
\overleftarrow{\mathrm{d}}\!\mu_{\mathrm s}(t)
&=
-\left[
-b\mu_{\mathrm s}(t)
+cX_t
-\frac{\sigma_Y^2}{R_{\mathrm f}(t)}
\left(\mu_{\mathrm f}(t)-\mu_{\mathrm s}(t)\right)
\right]\mathrm{d} t,
\\
\overleftarrow{\mathrm{d}}\!R_{\mathrm s}(t)
&=
-\left[
2\left(
-b+\frac{\sigma_Y^2}{R_{\mathrm f}(t)}
\right)R_{\mathrm s}(t)
-\sigma_Y^2
\right]\mathrm{d} t.
\end{aligned}
\label{eq:sec3_scalar_backward_smoother}
\end{equation}

\section{Derivation of Stationary Observable Statistics}
\label{app:stationary_statistics}

This appendix derives the analytical expressions for the stationary observable statistics used for the parameter choices in Section~\ref{sec:linear}. Throughout this appendix, we assume that the system has reached its stationary regime ($t\gg 0$ sufficiently large). Sections~\ref{app:two_dimensional_statistics} and~\ref{app:two_timescale_statistics} derive the stationary observable variance and decorrelation time for the two-dimensional linear model in Test I and the three-dimensional linear model with two hidden forcing timescales in Test II, respectively.

\subsection{Two-dimensional linear model (Test I)}
\label{app:two_dimensional_statistics}

For the stable linear model \eqref{eq:sec3_perfect_model}, write \(\mathbf{Z}_t=(X_t,Y_t)^\top\). Its drift matrix and noise covariance are
\[
\mathbf{A}:=\begin{pmatrix}-a&1\\c&-b\end{pmatrix},
\qquad \mathbf{Q}:=\operatorname{diag}(\sigma_X^2,\sigma_Y^2).
\]
Applying It\^o's formula to \(\mathbf{Z}_t\mathbf{Z}_t^\top\), taking expectations, and setting the resulting time derivative to zero shows that the stationary covariance \(\boldsymbol{\Sigma}=\mathbb E[\mathbf{Z}_t\mathbf{Z}_t^\top]\) satisfies the continuous-time Lyapunov equation~\cite{sarkka2019applied}
\begin{equation}
\mathbf{A}\boldsymbol{\Sigma}+\boldsymbol{\Sigma} \mathbf{A}^\top+\mathbf{Q}=\mathbf{0}.
\label{eq:appB_covariance_lyapunov}
\end{equation}
Writing \(V_X=\Sigma_{11}\) and equating the three distinct matrix entries gives
\begin{equation}
\begin{aligned}
-2aV_X+2\Sigma_{12}+\sigma_X^2&=0,\\
cV_X-(a+b)\Sigma_{12}+\Sigma_{22}&=0,\\
2c\Sigma_{12}-2b\Sigma_{22}+\sigma_Y^2&=0.
\end{aligned}
\label{eq:appB_covariance_balances}
\end{equation}
The first and third equations yield \(\Sigma_{12}=aV_X-\sigma_X^2/2\) and \(\Sigma_{22}=(c\Sigma_{12}+\sigma_Y^2/2)/b\). Substituting these into the second equation gives
\begin{equation}
V_X=\frac{\sigma_Y^2+[b(a+b)-c]\sigma_X^2}{2(a+b)(ab-c)},
\qquad
m_X:=\frac{\Sigma_{12}}{V_X}=a-\frac{\sigma_X^2}{2V_X}.
\label{eq:appB_covariance_solution}
\end{equation}
This recovers the stationary variance used in Section~\ref{sec:sec3_linear_systems}.

To obtain the decorrelation time, we need to compute the lagged covariance \(\operatorname{Cov}(\mathbf{Z}_{t+\tau},\mathbf{Z}_t)\). Applying the variation-of-constants formula~\cite{sarkka2019applied} to \eqref{eq:sec3_perfect_model} gives
\[
\mathbf{Z}_{t+\tau}=e^{\mathbf{A}\tau}\mathbf{Z}_t+
\int_t^{t+\tau}e^{\mathbf{A}(t+\tau-s)}\mathbf{Q}^{1/2}\,d\mathbf{W}_s.
\]
Here \(\mathbf{W}_s\) is the vector of independent Wiener processes. Since the stochastic integral has a deterministic integrand and involves only future Wiener increments, it is independent of \(\mathbf{Z}_t\) and has zero mean. Hence, for \(\tau\geq0\),
\begin{equation}
\begin{gathered}
\operatorname{Cov}(\mathbf{Z}_{t+\tau},\mathbf{Z}_t)=e^{\mathbf{A}\tau}\boldsymbol{\Sigma},\\
r(\tau):=\operatorname{Cov}(X_{t+\tau},X_t)
=\mathbf{e}_1^\top\operatorname{Cov}(\mathbf{Z}_{t+\tau},\mathbf{Z}_t)\mathbf{e}_1
=\mathbf{e}_1^\top e^{\mathbf{A}\tau}\boldsymbol{\Sigma} \mathbf{e}_1,
\end{gathered}
\label{eq:appB_lagged_covariance}
\end{equation}
where \(\mathbf{e}_1=(1,0)^\top\). The autocorrelation function of \(X\) is \(\rho_X(\tau)=r(\tau)/V_X\). With \(\mathbf{z}=\boldsymbol{\Sigma} \mathbf{e}_1\), squaring \(\rho_X(\tau)\) and integrating gives
\begin{equation}
\begin{aligned}
\mathcal T_X
&=\frac{1}{V_X^2}\int_0^\infty
\mathbf{z}^\top e^{\mathbf{A}^\top\tau}\mathbf{e}_1\mathbf{e}_1^\top e^{\mathbf{A}\tau}\mathbf{z}\,d\tau
=\frac{\mathbf{z}^\top \mathbf{P}\mathbf{z}}{V_X^2},\\
\mathbf{P}:&=\int_0^\infty e^{\mathbf{A}^\top\tau}\mathbf{e}_1\mathbf{e}_1^\top e^{\mathbf{A}\tau}\,d\tau.
\end{aligned}
\label{eq:appB_gramian_representation}
\end{equation}
Since \(\mathbf{A}\) has eigenvalues with negative real parts, \(e^{\mathbf{A}\tau}\to\mathbf{0}\) as \(\tau\to\infty\). Differentiating the integrand of the auxiliary matrix \(\mathbf{P}\) with respect to \(\tau\) and integrating over \([0,\infty)\) gives
\begin{equation}
\mathbf{A}^\top \mathbf{P}+\mathbf{P}\mathbf{A}
=\left[e^{\mathbf{A}^\top\tau}\mathbf{e}_1\mathbf{e}_1^\top e^{\mathbf{A}\tau}\right]_0^\infty
=-\mathbf{e}_1\mathbf{e}_1^\top.
\label{eq:appB_gramian_lyapunov}
\end{equation}
Writing the symmetric matrix \(\mathbf{P}=(p_{ij})\) with \(p_{12}=p_{21}\) and equating the \((2,2)\) and \((1,2)\) entries on both sides of \eqref{eq:appB_gramian_lyapunov} gives \(p_{12}=bp_{22}\) and \(p_{11}=[b(a+b)-c]p_{22}\). Substituting these relations into the \((1,1)\) equation determines \(p_{22}\), giving
\begin{equation}
\mathbf{P}=\frac{1}{2(a+b)(ab-c)}
\begin{pmatrix}
b(a+b)-c & b\\
b & 1
\end{pmatrix}.
\label{eq:appB_gramian_solution}
\end{equation}
Finally, using \(\mathbf{z}=\boldsymbol{\Sigma}\mathbf{e}_1=(V_X,\Sigma_{12})^\top\) and \(m_X=\Sigma_{12}/V_X\), we obtain \(\mathbf{z}/V_X=(1,m_X)^\top\). Substituting this into \eqref{eq:appB_gramian_representation} gives
\begin{equation}
\mathcal T_X
=p_{11}+2m_Xp_{12}+m_X^2p_{22}
=\frac{m_X^2+2bm_X+b(a+b)-c}{2(a+b)(ab-c)},
\label{eq:appB_squared_time_evaluation}
\end{equation}
which is the closed form in \eqref{eq:sec3_calibration_tcorr2}.

\subsection{Three-dimensional linear model (Test II)}
\label{app:two_timescale_statistics}

For the linear model with two forcing timescales \eqref{eq:sec3_multiscale_perfect}, we decompose the observed variable into three mutually uncorrelated components, \(X_t=X_t^{(0)}+X_t^{(1)}+X_t^{(2)}\). The observable equation can then be rewritten as:
\begin{equation}
\begin{aligned}
dX_t^{(0)}&=-aX_t^{(0)}\,dt+\sigma_X\,dW_t^X,\\
dX_t^{(i)}&=(-aX_t^{(i)}+Y_t^{(i)})\,dt,\qquad i=1,2.
\end{aligned}
\label{eq:appB_test2_component_dynamics}
\end{equation}
Here \(X_t^{(0)}\) is the component of \(X_t\) driven by the direct noise \(\sigma_X\,dW_t^X\), while \(X_t^{(i)}\) is the component driven by the hidden forcing \(Y_t^{(i)}\), for \(i=1,2\) (cf., \eqref{eq:sec3_X_memory}).

For \(i=0,1,2\), let \(V_i:=\operatorname{Var}(X_t^{(i)})\) denote the stationary variance of the component \(X_t^{(i)}\).
For \(i=1,2\), applying It\^o's formula \cite{sarkka2019applied} to \((Y_t^{(i)})^2\), \(X_t^{(i)}Y_t^{(i)}\), and \((X_t^{(i)})^2\), taking expectations, and using stationarity gives
\begin{equation}
\begin{aligned}
-2b_i\operatorname{Var}(Y_t^{(i)})+(\sigma_Y^{(i)})^2&=0,\\
-(a+b_i)\operatorname{Cov}(X_t^{(i)},Y_t^{(i)})+\operatorname{Var}(Y_t^{(i)})&=0,\\
-2aV_i+2\operatorname{Cov}(X_t^{(i)},Y_t^{(i)})&=0.
\end{aligned}
\label{eq:appB_test2_moment_balances}
\end{equation}
The second and third equations in \eqref{eq:appB_test2_moment_balances} give
\begin{equation}
\frac{\operatorname{Var}(Y_t^{(i)})}{a+b_i}=
\operatorname{Cov}(X_t^{(i)},Y_t^{(i)})
=aV_i.
\label{eq:appB_test2_covariance_variance_relation}
\end{equation}
The first equation in \eqref{eq:appB_test2_moment_balances} gives the hidden-state variance \(\operatorname{Var}(Y_t^{(i)})=(\sigma_Y^{(i)})^2/(2b_i)\). Substituting this into \eqref{eq:appB_test2_covariance_variance_relation} and dividing by \(a\) yields \(V_i=(\sigma_Y^{(i)})^2/[2ab_i(a+b_i)]\). For \(X_t^{(0)}\), applying It\^o's formula to \((X_t^{(0)})^2\) using \eqref{eq:appB_test2_component_dynamics}, taking expectations, and using stationarity gives
\begin{equation}
-2aV_0+\sigma_X^2=0.
\label{eq:appB_test2_direct_variance_balance}
\end{equation}
Hence, \(V_0=\sigma_X^2/(2a)\). Adding the variances of the three uncorrelated components gives the total stationary observable variance \(V_X=\operatorname{Var}(X_t)\):
\begin{equation}
\begin{aligned}
V_X=V_0+V_1+V_2,\qquad
V_0=\frac{\sigma_X^2}{2a},
\qquad
V_i=\frac{(\sigma_Y^{(i)})^2}{2ab_i(a+b_i)},\quad i=1,2.
\end{aligned}
\label{eq:appB_test2_variance}
\end{equation}
Here \(V_0\), \(V_1\), and \(V_2\), in the context of the study in Section \ref{sec:linear-test-ii}, are the variances of the components driven by the direct noise, slow hidden forcing, and fast hidden forcing, respectively.

To compute the decorrelation time, we first obtain the autocovariance of each component \(X_t^{(i)}\), \(i=0,1,2\), for \(\tau\geq0\). The component \(X_t^{(0)}\), driven by direct noise, is an OU process, with autocovariance
\[
r_0(\tau):=\operatorname{Cov}(X_{t+\tau}^{(0)},X_t^{(0)})=V_0e^{-a\tau}.
\]
For the components \(X_t^{(i)}\) driven by hidden forcing, \(i=1,2\), define \(r_i(\tau):=\operatorname{Cov}(X_{t+\tau}^{(i)},X_t^{(i)})\). Applying the variation-of-constants formula to the hidden-state equation in \eqref{eq:sec3_multiscale_perfect} gives
\[
Y_{t+\tau}^{(i)}=e^{-b_i\tau}Y_t^{(i)}
+\sigma_Y^{(i)}\int_t^{t+\tau}e^{-b_i(t+\tau-s)}\,dW_s^{(i)}.
\]
The stochastic integral has a deterministic integrand and involves only future Wiener increments, so it is independent of \(X_t^{(i)}\) and has zero mean. Taking the covariance of both sides with \(X_t^{(i)}\) and using \eqref{eq:appB_test2_covariance_variance_relation} therefore gives
\[
\operatorname{Cov}(Y_{t+\tau}^{(i)},X_t^{(i)})
=e^{-b_i\tau}\operatorname{Cov}(Y_t^{(i)},X_t^{(i)})
=aV_ie^{-b_i\tau}.
\]
Holding \(t\) fixed, differentiating \(r_i(\tau)\) with respect to \(\tau\) and substituting the component dynamics \eqref{eq:appB_test2_component_dynamics} gives
\[
\begin{aligned}
\frac{dr_i}{d\tau}
&=\operatorname{Cov}(-aX_{t+\tau}^{(i)}+Y_{t+\tau}^{(i)},X_t^{(i)})\\
&=-ar_i(\tau)+aV_ie^{-b_i\tau},
\qquad r_i(0)=V_i.
\end{aligned}
\]
Applying the variation-of-constants formula to this equation and evaluating the integral gives, for \(a\ne b_i\),
\begin{equation}
r_i(\tau)=V_i\frac{ae^{-b_i\tau}-b_ie^{-a\tau}}{a-b_i}.
\label{eq:appB_test2_component_autocovariance}
\end{equation}
The observable autocovariance, from properties of the covariance operator for uncorrelated sums, enjoys the following decomposition:
\begin{equation}
\begin{aligned}
r(\tau)&:=\operatorname{Cov}(X_{t+\tau},X_t)
=r_0(\tau)+r_1(\tau)+r_2(\tau)\\
&=C_ae^{-a\tau}+\sum_{i=1}^2C_ie^{-b_i\tau},
\qquad r(0)=V_X,
\end{aligned}
\label{eq:appB_test2_autocovariance}
\end{equation}
where
\[
C_i=\frac{aV_i}{a-b_i},
\qquad
C_a=V_0-\sum_{i=1}^2\frac{b_iV_i}{a-b_i}.
\]
The autocorrelation function is \(\rho_X(\tau)=r(\tau)/V_X\). Substituting \eqref{eq:appB_test2_autocovariance} into \(\mathcal T_X=\int_0^\infty\rho_X(\tau)^2\,d\tau\), expanding the square, and integrating each term using \(\int_0^\infty e^{-\lambda\tau}\,d\tau=1/\lambda\) for \(\lambda>0\) gives
\begin{equation}
\mathcal T_X=\frac{1}{V_X^2}
\left[
\frac{C_a^2}{2a}
+\sum_{i=1}^2\frac{C_i^2}{2b_i}
+2\sum_{i=1}^2\frac{C_aC_i}{a+b_i}
+\frac{2C_1C_2}{b_1+b_2}
\right].
\label{eq:appB_test2_Tcorr2}
\end{equation}
These expressions give the perfect-model statistics reported in Test II.

\bibliographystyle{aipnum4-1}
\bibliography{causal_fingerprintsNotes,references}

\end{document}